\PassOptionsToPackage{noend}{algpseudocode}
\PassOptionsToPackage{usenames,dvipsnames,table}{xcolor}
\documentclass[a4paper, 10pt, USenglish, reqno]{amsart}

\usepackage[%
	eqreset=section,
	thmreset=section,
	noload={caption,subcaption},
]{myPreamble}

\CreateNewTheorem[Fact]{fact}[dummythm]{Fact}
\crefformat{subsection}{#2\S #1#3}

	\let\proj\relax
	\DeclareMathOperator\proj{Proj}

	\newcommand{\C}{\mathcal C}
	\newcommand{\M}{\mathcal M}
	\renewcommand{\FBsmooth}{{\vphantom f\smash{\tilde f}}}%
	\renewcommand{\FBnonsmooth}{{\vphantom g\smash{\tilde g}}}%
	\renewcommand{\FBstepsize}{{\tilde\gamma}}%

	\makeatletter
		\newcommand\D{\@ifstar\@@D\@D}
		\newcommand\@D[1][h]{%
			\@ifnextchar_{\operatorname D}{%
				\ifstrempty{#1}{\operatorname D}{\operatorname D_{#1}}%
			}
		}
		\newcommand\@@D[1][\hat h]{%
			\@ifnextchar_{\operatorname D}{%
				\ifstrempty{#1}{\operatorname D}{\operatorname D_{#1}}%
			}
		}
		\newcommand{\h}{\@ifstar\@@h\@h}
		\newcommand{\@h}{\vphantom h\smash{\hat h}}
		\newcommand{\@@h}{\hat h}
		\newcommand{\Tkernel}{h}%
		\newcommand{\Tsmooth}{f}%
		\newcommand{\Tnonsmooth}{g}%
		\newcommand{\Tstepsize}{\gamma}%
		\newcommand{\FBE}{\@ifstar\@@FBE\@FBE}
		\newcommand{\@FBE}{\varphi_{\nicefrac h\gamma}^{f,g}}
		\newcommand{\@@FBE}{\varphi_{\nicefrac{1}{\tilde\gamma}}^{\tilde f,\tilde g}}
		\newcommand{\T}{\@ifstar\@@T\@T}
		\newcommand{\@T}{\operatorname T_{\nicefrac{\Tkernel}{\Tstepsize}}^{\Tsmooth,\Tnonsmooth}}
		\newcommand{\@@T}{\operatorname T_{\nicefrac{1}{\FBstepsize}}^{\FBsmooth,\FBnonsmooth}}
		\newcommand{\Res}{\operatorname R_{\nicefrac h\gamma}^{f,g}}
		\renewcommand{\Q}{Q_{\nicefrac h\gamma}^f}
	\makeatother

	\crefalias{AlgoLine}{line}%

	\let\oldCref\Cref
	\let\oldcref\cref
	\makeatletter
		\newcommand\algname{{\tt\textsc{Bella}}}
		\renewcommand\Cref[1]{%
			\IfStrEq{#1}{alg:Bella}{%
				\hyperref[alg:Bella]{\algname}
			}{%
				\oldCref{#1}%
			}%
		}%
		\renewcommand\cref[1]{%
			\IfStrEq{#1}{alg:Bella}{%
				\hyperref[alg:Bella]{\algname}%
			}{%
				\oldcref{#1}%
			}%
		}%
	\makeatother

\newif\ifshownotes\shownotestrue

\shownotesfalse
\showgrayfalse
\tikzexternaldisable

\newcommand\TheKeywords{%
	Nonsmooth nonconvex optimization,
	Bregman-Moreau envelope,
	Bregman forward-backward envelope,
	relative smoothness,
	KL inequality,
	nonlinear error bound,
	nonisolated local minima,
	superlinear convergence%
}
\newcommand\TheSubjclass{%
	90C06, 
	90C25, 
	90C26, 
	49J52, 
	49J53.
}
\newcommand\TheTitle{%
	A Bregman forward-backward linesearch algorithm for nonconvex composite optimization:\ifsiam\texorpdfstring{\\}{ }\else\ \fi superlinear convergence to nonisolated local minima%
}
\newcommand\TheShortTitle{%
	Bella: A superlinearly convergent Bregman forward-backward algorithm%
}
\newcommand\TheShortAuthor{%
	M. Ahookhosh, A. Themelis, and P. Patrinos%
}
\newcommand\TheFunding{%
	This work was supported by the \emph{Research Foundation Flanders (FWO)} research projects G086518N, G086318N, and G0A0920N;
	\emph{Research Council KU Leuven} C1 project No. C14/18/068;
	\emph{Fonds de la Recherche Scientifique --- FNRS and the Fonds Wetenschappelijk Onderzoek --- Vlaanderen} under EOS project no 30468160 (SeLMA)%
}

\ifsiam
	
	\headers{\TheShortTitle}{\TheShortAuthor}
	
	\title{%
		\TheTitle%
		\thanks{%
			Submitted to the editors \today.%
			\funding{\TheFunding}%
		}%
	}
	
	\author{%
		Masoud Ahookhosh\texorpdfstring{\footnotemark[2]}{}~
		Andreas Themelis\texorpdfstring{\footnotemark[2]}{}~ and
		Panagiotis Patrinos%
		\thanks{%
			\TheAddressKU\ {\tt
				\{%
					\href{mailto:masoud.ahookhosh@esat.kuleuven.be}{masoud.ahookhosh},%
					\href{mailto:andreas.themelis@esat.kuleuven.be}{andreas.themelis},%
					\href{mailto:panos.patrinos@esat.kuleuven.be}{panos.patrinos}%
				\}%
				\href{mailto:masoud.ahookhosh@esat.kuleuven.be,andreas.themelis@esat.kuleuven.be,panos.patrinos@esat.kuleuven.be}{@esat.kuleuven.be}%
			}%
		}%
	}%

\else

	\title[\TheShortTitle]{\sc\bfseries\TheTitle}
	\author[\TheShortAuthor]{%
		Masoud Ahookhosh,
		Andreas Themelis and
		Panagiotis Patrinos%
	}
	\thanks{\TheAddressKU\\\TheFunding}
	\subjclass{\TheSubjclass}
\fi

\begin{document}

	\ifams\else
		\maketitle
	\fi
	\begin{abstract}
		We introduce \Cref{alg:Bella}, a locally superlinearly convergent Bregman forward backward splitting method for minimizing the sum of two nonconvex functions, one of which satisfying a relative smoothness condition and the other one possibly nonsmooth.
		A key tool of our methodology is the Bregman forward-backward envelope (BFBE), an exact and continuous penalty function with favorable first- and second-order properties, and enjoying a nonlinear error bound when the objective function satisfies a {\L}ojasiewicz-type property.
		The proposed algorithm is of linesearch type over the BFBE along candidate update directions, and converges subsequentially to stationary points, globally under a KL condition, and owing to the given nonlinear error bound can attain superlinear convergence rates even when the limit point is a nonisolated minimum, provided the directions are suitably selected.
	\end{abstract}
	
	\ifams
		\maketitle
	\else
	\fi


	\section{Introduction}\label{sec:Introduction}
		In this paper, we address the composite minimization problem
\begin{equation}\label{eq:Pinit}
	\minimize \varphi(x)\equiv f(x)+g(x)
\quad\stt{}
	x\in\overline C.
\end{equation}
Here, \(C\) is an open convex set (\(\overline C\) denotes the closure of \(C\)), \(g\) is proper and lower semicontinuous (lsc), and \(f\) is relatively smooth with respect to a Legendre kernel \(h\) (see \cref{sec:Smoothness}) with \(\dom\nabla h=C\)
(for detailed assumptions we refer to \cref{sec:FBS}).
Despite its simple structure, \eqref{eq:Pinit} encompasses a variety of optimization problems frequently encountered in scientific areas such as signal and image processing, machine learning, and inverse problems \cite{bauschke2016descent,bolte2018first,hanzely2018accelerated,lu2018relatively,mukkamala2019beyond}. 
The notion of \DEF{Lipschitz-like convexity} was recently discovered in the seminal \cite{bauschke2016descent} as a generalization of the Lipschitz smoothness condition, which is latter named as \DEF{relative smoothness} in \cite{lu2018relatively}.
Studying optimization problems involving relatively smooth functions has received much attention in the last few years \cite{bauschke2016descent,bolte2018first,dragomir2019fast,hanzely2018fastest,hanzely2018accelerated,lu2018relatively,nesterov2018implementable,teboulle2018simplified}.
In our setting \eqref{eq:Pinit}, since $f$ is relatively smooth and $g$ is nonsmooth nonconvex, we can cover a wide spectrum of  applications. In the Euclidean setting, there are plenty of optimization algorithms that can handle composite minimization of the form \eqref{eq:Pinit}, such as \cite{ahookhosh2019accelerated,ahookhosh2018optimal,beck2009fast,bolte2017error,nesterov2013gradient,themelis2019acceleration} for convex problems and \cite{attouch2010proximal,bolte2014proximal,bot2016inertial,bot2016inertialTseng,
bot2019proximal,fukushima1981generalized,stella2017forward,themelis2018forward} for nonconvex problems.

One of the most significant discussions in the field of numerical optimization has been related to designing iterative schemes guaranteeing a superlinear convergence rate; see, e.g., \cite{nocedal2006numerical} for many algorithms attaining a superlinear convergence rate for smooth problems and \cite{fukushima1981generalized,fukushima1996globally,themelis2018forward} for other related works in the nonconvex nonsmooth setting.
In most of these attempts, the key element is the so-called Dennis-Moré condition \cite{dennis1974characterization,dennis1977quasi} which guarantees superlinear convergence to an isolated critical point of the objective function.
However, there are many applications that have nonisolated critical points such as low-rank matrix completion \cite{tanner2016low}, low-rank matrix recovery \cite{bhojanapalli2016global}, phase retrieval \cite{sun2018geometric}, and deep learning \cite{kawaguchi2016deep}.
Up to now, besides some attempts for minimizing smooth nonlinear least-squares problems (see, e.g., \cite{ahookhosh2019local,ahookhosh2020finding,izmailov2014newton} and references therein) far too little attention has been paid to the superlinear convergence to nonisolated critical points for nonconvex nonsmooth problems.

		\subsection{Related work}\label{sec:stateArt}
			In order to guarantee convergence, most first-order methods for problem \eqref{eq:Pinit} in the nonconvex setting require Lipschitz differentiability of \(f\); however, there are plenty of examples that fail to comply with this assumption; see, e.g., \cite{ahookhosh2020block,ahookhosh2019multi,bauschke2016descent,dragomir2019quartic,lu2018relatively,mukkamala2019beyond}.
Recently, from the seminal work conducted in \cite{bauschke2016descent} it emerged that the Lipschitz smoothness assumption of \(f\) can be relaxed by introducing the notion of \DEF{relative smoothness} (see \cref{def:relSmooth}), which was further developed in \cite{lu2018relatively,van2017forward}.
Assuming the convexity of \(f\) and \(g\) and relative smoothness of \(f\), a Bregman proximal gradient method was proposed in \cite{bauschke2016descent}, while primal and dual algorithms were developed in \cite{lu2018relatively}.
More recently, in the convex setting, \cite{nesterov2018implementable} proposed an accelerated tensor method, and \cite{hanzely2018fastest,hanzely2018accelerated} suggested a Nesterov-type accelerated method and a stochastic mirror descent method.

The developments of relative smoothness also led to a renewed interest in theory and algorithms of nonconvex optimization.
Recently, \cite{bolte2018first} extended the results of \cite{bauschke2016descent} for the Bregman proximal gradient method, and \cite{teboulle2018simplified} discussed several first-order algorithms.
More recently, the linear convergence of the gradient method for relatively smooth functions was studied in \cite{bauschke2019linear}.
In \cite{ochs2019nonsmooth}, a generic Bregman linesearch was proposed for not necessarily Lipschitz smooth problems, which covers the Bregman forward-backward splitting as a special case for nonconvex smooth \(f\) and convex \(g\).
In \cite{mukkamala2019convex}, some Bregman proximal gradient algorithms with inertial effects were presented in the nonconvex setting for smooth \(f\) and hypoconvex \(g\).
Furthermore, a stochastic convex model-based minimization algorithm was proposed in \cite{davis2018stochastic} for hypoconvex functions under relative smoothness and high-order growth conditions.
The notion of relative smoothness was further extended to its block version \cite{ahookhosh2020block,ahookhosh2019multi} opening the possibility of alternating minimization algorithms in the fully nonconvex setting.
To the best of our knowledge, apart from the latter papers, there have not been many attempts to deal with \eqref{eq:Pinit} in the relative smooth and fully nonconvex setting.

 		\subsection{Contribution}
			Our main goal is to design an efficient framework for addressing the structured nonconvex problem \eqref{eq:P} with superlinear guarantees of convergence, even when the limit point is a nonisolated local minimum.
We aim at devising a linesearch strategy that globalizes the convergence of fast local methods, stemming for instance from Newton-type schemes.
The lack of differentiability in problem \eqref{eq:Pinit} makes classical smooth optimization methodologies, such as Armijo backtracking, not applicable.
Nevertheless, favorable properties of the here introduced Bregman forward-backward envelope lead to the \emph{Bregman EnveLope Linesearch Algorithm} (\cref{alg:Bella}), which overcomes said limitations exclusively by means of Bregman forward-backward operations.
Our contribution can be summarized as follows.
\begin{enumerate}[wide, labelwidth=!, labelindent=0pt]
\item{\it Bregman-Moreau envelope analysis.}
	We provide new insights on the Bregman-Moreau envelope complementing the ones in \cite{kan2012moreau,bauschke2018regularizing,laude2020bregman}.
	Among these, we highlight properties of fixed points, a local equivalence with the forward-backward envelope \cite{patrinos2013proximal}, and ultimately a second-order differentiability result with classical generalized differentiability tools \cite{poliquin1995second,rockafellar2011variational}.
\item{\it Bregman forward-backward splitting and linesearch extension.}
	We highlight a connection between Bregman forward-backward mapping and Bregman proximal mapping (\Cref{thm:Legendre}), revealing that the proximal point algorithm is as general as the forward-backward splitting in the Bregman setting.
	Correspondingly, we introduce a Bregman forward-back\-ward envelope function (BFBE) which is at the core of \cref{alg:Bella}, a linesearch algorithm that can globalize the convergence of fast local methods for problem \eqref{eq:Pinit}.
\item{\it Superlinear convergence to nonisolated critical points.}
	Differentiability properties of the BFBE (\cref{thm:1st}) allow us to link the Kurdyka-{\L}ojasiewicz (KL) property to a nonlinear error bound involving the distance to sublevel sets (\Cref{thm:EB}).
	This observation is the key for showing that \cref{alg:Bella} can converge superlinearly to local minima even if nonisolated (\Cref{thm:tau1}), when the directions are suitably selected.
	More generally, global and linear convergence are shown by means of the KL-inequality (\cref{thm:global,thm:linear}).
\end{enumerate}

		\subsection{Paper organization}
			The rest of the paper is organized as follows.
\Cref{sec:Preliminaries} introduces the used notation and some known facts.
In \Cref{sec:prox}, we investigate some properties of the Bregman proximal mapping and the Bregman-Moreau envelope, which we then use in \Cref{sec:FBS} to derive similar results for the Bregman forward-backward mapping and the corresponding \emph{envelope} function, the Bregman forward-backward envelope (BFBE), which is a key tool of our analysis.
In \Cref{sec:Algorithm}, we introduce \Cref{alg:Bella}, a linesearch algorithm on the BFBE and show its convergence properties.
\Cref{sec:Conclusion} concludes the paper.

	\section{Preliminaries}\label{sec:Preliminaries}

		\subsection{Notation}
			The extended-real line is denoted by \(\Rinf\coloneqq\R\cup\set\infty\).
The open and closed balls of radius \(r\geq 0\) centered at \(x\in\R^n\) are denoted as \(\ball xr\) and \(\cball xr\), respectively.
We say that \(\seq{x^k}\subset\R^n\) converges at \(R\)-linear rate (to a point \(x_\star\)) if there exists \(c>0\) and \(\rho\in(0,1)\) such that \(\|x^k-x_\star\|\leq c\rho^k\) holds for every \(k\).
The distance of a point \(x\in\R^n\) to a nonempty set \(S\subseteq\R^n\) is given by \(\dist(x,S)=\inf_{z\in S} \|z-x\|\).
The closure and interior of \(S\) are respectively denoted as \(\interior S\) and \(\overline S\).

A function \(\func{f}{\R^n}{\Rinf}\) is \DEF{proper} if \(f\not\equiv\infty\), in which case its \DEF{domain} is defined as the set \(\dom f\coloneqq\set{x\in\R^n}[f(x)<\infty]\).
For \(\alpha\in\R\),
\(
	[f\leq\alpha]
{}\coloneqq{}
	\set{x\in\R^n}[f(x)\leq\alpha]
\)
is the \DEF{\(\alpha\)-(sub)level set} of \(f\); \([f\geq\alpha]\), \([f=\alpha]\), etc. are defined accordingly.
We say that \(f\) is \DEF{level bounded} if \([f\leq\alpha]\) is bounded for all \(\alpha\in\R\).
A point \(x_\star\in\dom f\) is a \DEF{local minimum} for \(f\) if \(f(x)\geq f(x_\star)\) holds for all \(x\) in a neighborhood of \(x_\star\).
If the inequality can be strengthened to \(f(x)\geq f(x_\star)+\tfrac\mu2\|x-x_\star\|^2\) for some \(\mu>0\), then \(x_\star\) is a \DEF{strong local minimum}.
The \DEF{convex conjugate} of \(f\) is denoted as \(\conj f\coloneqq\sup_z\set{\innprod{{}\cdot{}}{z}-h(z)}\).

A vector \(v\in\partial f(x)\) is a \DEF{subgradient} of \(f\) at \(x\), where \(\partial f(x)\) is the \DEF{(limiting) subdifferential}
\[
	\partial f(x)
{}\coloneqq{}
	\set{v\in\R^n}[
		\exists\seq{x^k,v^k}~\text{s.t.}~ x^k\to x,~f(x^k)\to h(x),~
		\hat\partial f(x^k)\ni v^k\to v
	],
\]
and \(\hat\partial f(x)\) is the set of \DEF{regular subgradients} of \(f\) at \(x\), namely vectors \(v\in\R^n\) such that
\(
	\liminf_{\limsubstack{z&\to&x\\z&\neq&x}}{
		\frac{
			f(z)-f(x)-\innprod{v}{z-x}
		}{
			\|z-x\|
		}
	}
{}\geq{}
	0
\).
Following the terminology of \cite{rockafellar2011variational}, we say that \(\func{f}{\R^n}{\Rinf}\) is \DEF{strictly continuous at \(\bar x\)} if
\(
	\lip f(\bar x)
{}\coloneqq{}
	\limsup_{
		\substack{y,z\to\bar x\\y\neq z}
	}{
		\frac{|f(y)-f(z)|}{\|y-z\|}
	}
{}<{}
	\infty
\),
and \DEF{strictly differentiable at \(\bar x\)} if \(\nabla f(\bar x)\) exists and satisfies
\(
	\lim_{
		\substack{y,z\to\bar x\\y\neq z}
	}{
		\frac{f(y)-f(z) - \innprod{\nabla f(\bar x)}{y-z}}{\|y-z\|}
	}
{}={}
	0
\).
If \(f\) is everywhere strictly continuous on an open set \(\mathcal U\), then its gradient exists almost everywhere on \(\mathcal U\), and as such its \DEF{Bouligand subdifferential}
\[
	\partial_Bf(x)
{}\coloneqq{}
	\set{v}[
		\exists x^k\to x \text{ with } \nabla f(x^k)\to v
	]
\]
is nonempty and compact for all \(x\in\mathcal U\) \cite[Thm. 9.61]{rockafellar2011variational}.
\(\C^k(\mathcal U)\) is the set of functions \(\mathcal U\to\R\) which are \(k\) times continuously differentiable.
We simply write \(\C^k\) if \(\mathcal U\) is clear from context.
The notation \(\ffunc{T}{\R^n}{\R^n}\) indicates a point-to-set mapping, whose domain and range are respectively defined as
\(
	\dom T=\set{x\in\R^n}[T(x)\neq\emptyset]
\)
and
\(
	\range T=\bigcup_{x\in\R^n}T(x)
\).

		\subsection{Relative smoothness}\label{sec:Smoothness}
			Here, after giving some definitions, we establish necessary facts regarding relative smoothness.

\begin{defin}%
	Let \(\func{h}{\R^n}{\Rinf}\) be a proper, lsc, convex function with \(\interior\dom h\neq\emptyset\) and such that  \(h\in\C^1(\interior\dom h)\).
	Then, \(h\) is said to be
	\begin{enumerate}
	\item
		a \DEF{kernel function} if it is \DEF{\(1\)-coercive}, \ie \(\lim_{\|x\|\to\infty}\nicefrac{h(x)}{\|x\|}=\infty\);
	\item
		\DEF{essentially smooth}, if \(\|\nabla h(x_k)\|\to\infty\) for every sequence \(\seq{x_k}\subseteq\interior\dom h\) converging to a boundary point of \(\dom h\);
	\item
		of \DEF{Legendre} type if it is essentially smooth and strictly convex.
	\end{enumerate}
\end{defin}

\begin{defin}[Bregman distance \cite{bregman1967relaxation}]%
	For a convex function \(\func{h}{\R^n}{\Rinf}\) continuously differentiable on \(\interior\dom h\neq\emptyset\), the \DEF{Bregman distance} \(\func{\D}{\R^n\times\R^n}{\Rinf}\) is given by
	\begin{equation}\label{eq:bregman}
		\D(x,y)
	{}\coloneqq{}
		\begin{ifcases}
			h(x)-h(y)-\innprod{\nabla h(y)}{x-y} & y\in\interior\dom h
		\\
			\infty\otherwise.
		\end{ifcases}
	\end{equation}
\end{defin}

If \(h\) is a strictly convex kernel function, then \(\D\) serves as a pseudo-distance, having \(\D\geq 0\) and \(\D(x,y)=0\) iff \(x=y\in\interior\dom h\).
In general, however, \(\D\) is nonsymmetric and fails to satisfy the triangular inequality.
There are many popular kernel functions such as energy, Boltzmann-Shannon entropy, Fermi-Dirac entropy and so on leading to various Bregman distances that appear in many applications; see, e.g., \cite[Ex. 2.3]{bauschke2018regularizing}.

\begin{rem}
	The following assertions hold: 
	\begin{enumerate}
	\item\label{thm:inv1}%
		If \(\func{h}{\R^n}{\Rinf}\) is a Legendre kernel, then \(\conj h\in\C^1(\R^n)\) is strictly convex.
		In fact, \(\func{\nabla h}{\interior\dom h}{\R^n}\) is a bijection with \(\nabla h^{-1}=\nabla\conj h\)
		\cite[Thm. 26.5, Cor. 13.3.1]{rockafellar1970convex}.
	\item\label{thm:conjC2}%
		If \(h\in\C^2\) is of Legendre type and \(\nabla^2h\succ0\) on \(\interior\dom h\), then \(\conj h\in\C^2(\R^n)\)
		\cite{rockafellar1977higher}.
	\item\label{thm:Dhcoercive}%
		\(\D({}\cdot{},x)\) and \(\D(x,{}\cdot{})\) are level bounded locally uniformly in \(x\) on \(\interior\dom h\times\interior\dom h\)
		\cite[Lem. 7.3(v)-(viii)]{bauschke2001essential}.\footnote{%
			Although \cite{bauschke2001essential} only states level boundedness, a trivial modification of the proof shows local uniformity too.%
		}%
	\end{enumerate}
	Moreover, for any open convex set \(\mathcal U\subseteq\interior\dom h\) the following hold:
	\begin{enumerate}
	\item\label{thm:hstrcvx}%
		If \(\nabla h\) is \(\tilde\sigma_h\)-strongly monotone on \(\mathcal U\), then
		\(
			\D(y,x)
		{}\geq{}
			\frac{\tilde\sigma_h}{2}\|y-x\|^2
		\)
		for all \(x,y\in\mathcal U\).
	\item\label{thm:hC11}%
		If \(\nabla h\) is \(\tilde L_h\)-Lipschitz on \(\mathcal U\), then
		\(
			\D(y,x)
		{}\leq{}
			\frac{\tilde L_h}{2}\|y-x\|^2
		\)
		for all \(x,y\in\mathcal U\).
	\hfill\qedsymbol
	\end{enumerate}\let\qedsymbol\relax
\end{rem}

We will sometimes require properties such as Lipschitz differentiability or strong convexity to hold locally, where locality amounts to the existence for any point of a convex neighborhood in which such property holds.

\begin{defin}[relative smoothness {\cite{bauschke2016descent}}]%
\label{def:relSmooth}
	We say that a proper, lsc function \(\func{f}{\R^n}{\Rinf}\) is \DEF{smooth relative to} a kernel \(\func{h}{\R^n}{\Rinf}\) if \(\dom f\supseteq\dom h\), and there exists \(L_f\geq0\) such that \(L_fh\pm f\) are convex functions on \(\interior\dom h\).
	We will simply say that \(f\) is relatively smooth when \(h\) is clear from context, or \(L_f\)-relatively smooth to make the modulus \(L_f\) explicit.
\end{defin}

\begin{prop}%
	Let \(f\) be smooth relative to a kernel \(h\).
	Then, the following hold:
	\begin{enumerate}
	\item\label{thm:fC1}
		\(f\in\C^1(\interior\dom h)\);
	\item\label{thm:fC11}
		if \(h\) is Lipschitz differentiable on an open set \(\mathcal U\), then so is \(f\).
	\end{enumerate}
	\begin{proof}
		\begin{proofitemize}
		\item\ref{thm:fC1}~
			Convexity of \(L_fh\pm f\) and continuous differentiability of \(h\) on \(\interior\dom h\) ensure through \cite[Ex. 8.20(b) and Cor. 9.21]{rockafellar2011variational} that both \(f\) and \(-f\) are subdifferentially regular on \(\interior\dom h\), in the sense of \cite[Def. 7.25]{rockafellar2011variational}, with \(\hat\partial f\) and \(\hat\partial(-f)\) both nonempty.
			The proof now follows from \cite[Thm. 9.18(d)]{rockafellar2011variational}.
		\item\ref{thm:fC11}~
			Let \(\tilde L_h\) be a Lipschitz modulus for \(\nabla h\) on \(\mathcal U\).
			Convexity of \(L_fh+f\) yields
			\begin{align*}
				\innprod{
					\nabla f(x)-\nabla f(y)
				}{
					x-y
				}
			{}\geq{} &
				-L_f\innprod{
					\nabla h(x)-\nabla h(y)
				}{
					x-y
				}
			{}\geq{}
				-L_f\tilde L_h\|x-y\|^2
			\shortintertext{%
				for \(x,y\in\mathcal U\), while due to concavity of \(f-L_fh\) it holds that
			}
				\innprod{
					\nabla f(x)-\nabla f(y)
				}{
					x-y
				}
			{}\leq{} &
				L_f\innprod{
					\nabla h(x)-\nabla h(y)
				}{
					x-y
				}
			{}\leq{}
				L_f\tilde L_h\|x-y\|^2.
			\end{align*}
			The two inequalities together prove that \(\nabla f\) is \(\tilde L_f\)-Lipschitz on \(\mathcal U\) with \(\tilde L_f=L_f\tilde L_h\).
			\qedhere
		\end{proofitemize}
	\end{proof}
\end{prop}

The proof of the following result is a a simple adaptation of that of \cite[Prop. 1.1]{lu2018relatively}.
\begin{prop}[characterization of relative smoothness]\label{thm:sigmaL}%
	The following are equivalent for a proper lsc function \(\func f{\R^n}{\Rinf}\):
	\begin{enumerateq}
	\item
		\(f\) is \(L_f\)-smooth relative to \(h\);
	\item\label{thm:LipIneq}%
		\(
			\bigl|f(y)-f(x)-\innprod{\nabla f(x)}{y-x}\bigr|
		{}\leq{}
			L_f\D_h(y,x)
		\)
		for all \(x,y\in\interior\dom h\);
	\item
		\(
			\bigl|\innprod{\nabla f(x)-\nabla f(y)}{x-y}\bigr|
		{}\leq{}
			L_f\innprod{\nabla h(x)-\nabla h(y)}{x-y}
		\)
		for all \(x,y\in\interior\dom h\);
	\item
		\(
			-L_f\nabla^2h
		{}\preceq{}
			\nabla^2 f
		{}\preceq{}
			L_f\nabla^2h
		\)
		on \(\interior\dom h\), provided that \(f,h\in\C^2(\interior\dom h)\).
	\end{enumerateq}
\end{prop}
%
%

	\section{Bregman proximal mapping and Moreau envelope}\label{sec:prox}
		This section is devoted to the analysis of the Bregman proximal mapping and the Bregman-Moreau envelope which we will need in the coming sections.
Our main results include their local equivalence with Euclidean objects, namely the forward-backward mapping and the forward-backward envelope (\cref{thm:Moreau=FBE}), respectively, which in turn will be used to derive novel second-order properties (\Cref{thm:Moreau:2nd}) as simple byproducts of similar results available in the literature.

Given a Legendre function \(\func{\h}{\R^n}{\Rinf}\), the \DEF{Bregman proximal mapping} of \(\func{\varphi}{\R^n}{\Rinf}\) with kernel \(\h\) and stepsize \(\gamma>0\) is the set-valued map \(\ffunc{\prox_{\gamma\varphi}^{\h}}{\interior\dom\h}{\R^n}\) given by
\begin{align}\label{eq:bprox}
	\prox_{\gamma\varphi}^{\h}(x)
{}={}
	\prox_{\varphi}^{\nicefrac{\h}{\gamma}}(x)
{}\coloneqq{} &
	\argmin_{z\in\R^n}\set{
		\varphi(z) + \tfrac{1}{\gamma}\D*(z,x)
	},
\shortintertext{
	and the corresponding \DEF{Bregman-Moreau envelope} is \(\func{\varphi^{\nicefrac{\h}{\gamma}}}{\R^n}{[-\infty,\infty]}\) defined as
}
\label{eq:bdre}
	\varphi^{\nicefrac{\h}{\gamma}}(x)
{}\coloneqq{} &
	\inf_{z\in\R^n}\set{\varphi(z) + \tfrac{1}{\gamma}\D*(z,x)}.
\end{align}

The first equality in \eqref{eq:bprox} owes to the invariance of the \(\argmin\) operator under positive scalings, whereas the notation in \eqref{eq:bdre} is justified from the identity \(\frac1\gamma\D*=\D*_{\nicefrac{\h}{\gamma}}\).
Although writing \(\prox_{\varphi}^{\nicefrac{\h}{\gamma}}\) better reflects the kinship with the envelope \(\varphi^{\nicefrac{\h}{\gamma}}\), the (equally well-posed) notation \(\prox_{\gamma\varphi}^{\h}\) is more consistent with the classical \(\prox_{\gamma\varphi}\) widely adopted in the Euclidean setting, that is, when \(\h=\tfrac12\|{}\cdot{}\|^2\).
For the sake of maintaining this consistency and for simplicity of exposition, the kernel \(\h\) will be omitted in the Euclidean case; we will thus write \(\prox_{\gamma\varphi}\) and \(\varphi^{\nicefrac1\gamma}\) to indicate the Euclidean proximal map and Moreau envelope of \(\varphi\) with stepsize \(\gamma\).

The proximal mapping is a fundamental building block of splitting algorithms, and the corresponding envelope function offers an extremely valuable tool for the convergence analysis of such schemes.
For this reason, a lot of research has been devoted to the study of their properties in the Euclidean case.
It is well known, for instance, that when \(\varphi\) is a (proper, lsc, and) convex function, its (Euclidean) proximal mapping is (single-valued and) nonexpansive and its Moreau envelope is convex and Lipschitz differentiable.
More generally, even in the nonconvex setting, the Moreau envelope is a continuous function sharing infimum and minimizers with \(\varphi\), and local smoothness properties have been established with variational analysis tools such as prox-regularity and epigraphical differentiation.

As detailed in the \nameref{sec:Introduction}, the extension to non-Euclidean kernels \(\h\) offers a significant additional degree of flexibility.
Furthermore, the Bregman proximal mapping can encapsulate an entire splitting algorithm, as is the case of the (Euclidean) proximal-gradient operator \(\prox_{\gamma g}(\id-\gamma\nabla f)\) that can equivalently be expressed as \(\prox_{f+g}^{\nicefrac{\h}{\gamma}}\), where \(\h\coloneqq\tfrac{1}{2\gamma}\|{}\cdot{}\|^2-f\).
In fact, it will be shown in \cref{thm:Legendre} that this is still true even for the proximal-gradient operator with an arbitrary Bregman kernel.
In other words, from a theoretical standpoint \emph{the proximal-gradient scheme offers no advantage in generality over the proximal point algorithm in the Bregman setting}.

This awareness emphasizes the importance of studying the proximal mapping in full generality.
Nonetheless, in the nonconvex setting, there is a big discrepancy between the well-studied Euclidean setting and the less mature Bregman generalization.
In an attempt to partially fill this gap, this section complements the analysis of \cite{kan2012moreau,laude2020bregman} for the proximal mapping and the Moreau envelope in the Bregman setting.
The extension --- or better, the \emph{``translation''} --- of all the results for the proximal gradient will then be derived as simple byproducts in the following section.
As the \emph{``translation''} entails a change of Bregman metric, in order to avoid confusion we use the hat version \(\h\) in this section and reserve the notation \(h\) for the Legendre kernel of the proximal-gradient mapping.
We begin by introducing the notion of Bregman-type prox-boundedness, which is a technical requirement ensuring the well definedness of the proximal map and the properness of the Moreau envelope.

\begin{defin}[\(\h\)-prox-boundedness]%
	Given a kernel \(\h\), a function \(\func{\varphi}{\R^n}{\Rinf}\) is said to be \DEF{\(\h\)-prox-bounded} if there exists \(\gamma>0\) such that \(\varphi^{\nicefrac{\h}{\gamma}}(x)>-\infty\) for some \(x\in\R^n\).
	The supremum of the set of all such \(\gamma\) is the threshold \(\gamma_\varphi^{\h}\) of the \(\h\)-prox-boundedness, \ie
	\[
		\gamma_\varphi^{\h}
	{}\coloneqq{}
		\sup\set{\gamma>0}[\exists x\in\R^n~\text{s.t.}~\varphi^{\nicefrac{\h}{\gamma}}(x)>-\infty].
	\]
\end{defin}

Note that whenever a proper and lsc function \(\varphi\) is lower bounded by an affine function on \(\dom h\) (as is the case when \(\varphi\) is convex or lower bounded on \(\dom h\)) then it is \(\h\)-prox-bounded with threshold \(\gamma_\varphi^{\h}=\infty\).
For more general functions, instead, the threshold plays a central role in dictating the range of feasible stepsizes \(\gamma\).

\begin{fact}[{regularity properties of \(\prox_{\gamma\varphi}^{\h}\) and \(\varphi^{\nicefrac{\h}{\gamma}}\) \cite[Thm.s 2.2, 2.3 and 2.4]{kan2012moreau}}]\label{thm:facts}
	Let \(\h\) be a Legendre kernel and \(\func{\varphi}{\R^n}{\Rinf}\) proper and lsc.
	For every \(\gamma\in(0,\gamma_\varphi^{\h})\) the following hold:
	\begin{enumerate}
	\item\label{thm:Moreau:dom}%
		\(\dom\varphi^{\nicefrac{\h}{\gamma}}=\dom\prox_{\gamma\varphi}^{\h}=\interior\dom\h\);
	\item\label{thm:prox:range}
		\(\range\prox_{\gamma\varphi}^{\h}\subseteq\dom\varphi\cap\dom\h\);
	\item\label{thm:proxosc}%
		\(\prox_{\gamma\varphi}^{\h}\) is locally bounded, compact-valued, and osc in \(\interior\dom\h\);
	\item\label{thm:MoreauC0}%
		\(\varphi^{\nicefrac{\h}{\gamma}}\) is real valued and continuous on \(\interior\dom\h\); in fact, it is locally Lipschitz if so is \(\nabla\h\).
	\end{enumerate}
\end{fact}

Next, we furnish some elementary connections between \(\varphi\) and its envelope \(\varphi^{\nicefrac{\h}{\gamma}}\).

\begin{prop}[relation between \(\varphi\) and \(\varphi^{\nicefrac{\h}{\gamma}}\)]\label{thm:Moreau:basic}%
	Let \(\h\) be a Legendre kernel and \(\func{\varphi}{\R^n}{\Rinf}\) proper and lsc.
	Then, for every \(\gamma\in(0,\gamma_\varphi^{\h})\)
	\begin{enumerate}
	\item\label{thm:Moreau:eq}
		\(
			\varphi(\bar x)
			{}+{}
			\frac1\gamma\D*(\bar x,x)
		{}={}
			\varphi^{\nicefrac{\h}{\gamma}}(x)
		{}\leq{}
			\varphi(x)
		\)
		for \(x\in\interior\dom\h\) and \(\bar x\in\prox_{\gamma\varphi}^{\h}(x)\), with
		\(
			\varphi^{\nicefrac{\h}{\gamma}}(x)
		{}={}
			\varphi(x)
		\)
		iff \(x\in\prox_{\gamma\varphi}^{\h}(x)\).
	\end{enumerate}
	Moreover, if \(\range\prox_{\gamma\varphi}^{\h}\subseteq\interior\dom\h\), then the following also hold:
	\begin{enumerate}[resume]
	\item\label{thm:Moreau:inf}
		\(
			\inf\varphi^{\nicefrac{\h}{\gamma}}
		{}={}
			\inf_{\interior\dom\h}\varphi
		\)
		and
		\(
			\argmin\varphi^{\nicefrac{\h}{\gamma}}
		{}={}
			\argmin_{\interior\dom\h}\varphi
		\);
	\item\label{thm:Moreau:lb}
		\(\varphi^{\nicefrac{\h}{\gamma}}\) is level bounded iff \(\varphi\) is level bounded on \(\interior\dom\h\).
	\end{enumerate}
	\begin{proof}
		\begin{proofitemize}
		\item\ref{thm:Moreau:eq}~
			The first equality is the definition of the Bregman proximal map and Moreau envelope; the inequality follows by considering \(z=x\) in the subproblem \eqref{eq:bprox} defining \(\varphi^{\nicefrac{\h}{\gamma}}\).
			In turn, the ``iff'' condition owes to the fact that \(\D*(z,x)=0\) iff \(z=x\) for all \(x,z\in\interior\dom\h\).
		\item\ref{thm:Moreau:inf}~
			It follows from assertion \ref{thm:Moreau:eq} that \(\inf\varphi^{\nicefrac{\h}{\gamma}}\leq\inf_{\interior\dom\h}\varphi\).
			Let a sequence \(\seq{x^k}\) be such that \(\varphi^{\nicefrac{\h}{\gamma}}(x^k)\to\inf\varphi^{\nicefrac{\h}{\gamma}}\) as \(k\to\infty\).
			Then, taking \(\bar x^k\in\prox_{\gamma\varphi}^{\h}(x^k)\subseteq\interior\dom\h\), assertion \ref{thm:Moreau:eq} ensures that
			\(
				\liminf_{k\to\infty}\varphi(\bar x^k)
			{}\leq{}
				\inf_{\interior\dom\h}\varphi
			\),
			hence the claimed equivalence of infima.
			If \(x\in\argmin\varphi^{\nicefrac{\h}{\gamma}}\), necessarily \(x\in\dom\varphi^{\nicefrac{\h}{\gamma}}=\interior\dom\h\) and there thus exists \(\bar x\in\prox_{\gamma\varphi}^{\h}(x)\), (cf. \cref{thm:facts}), which satisfies \(\bar x\in\interior\dom\h\) by assumption.
			Then,
			\[
				\varphi(\bar x)
			{}\leq{}
				\varphi^{\nicefrac{\h}{\gamma}}(x)
				{}-{}
				\tfrac1\gamma\D*(\bar x,x)
			{}={}
				\inf\varphi^{\nicefrac{\h}{\gamma}}
				{}-{}
				\tfrac1\gamma\D*(\bar x,x)
			{}={}
				\inf_{\interior\dom\h}\varphi
				{}-{}
				\tfrac1\gamma\D*(\bar x,x),
			\]
			where the inequality follows from assertion \ref{thm:Moreau:eq}.
			Therefore \(x=\bar x\in\argmin_{\interior\dom\h}\varphi\).
			Similarly, for \(x\in\argmin_{\interior\dom\h}\varphi\) it follows from assertion \ref{thm:Moreau:eq} that
			\(
				\varphi^{\nicefrac{\h}{\gamma}}(x)
			{}\leq{}
				\varphi(x)
			{}={}
				\inf_{\interior\dom\h}\varphi
			{}={}
				\inf\varphi^{\nicefrac{\h}{\gamma}}
			\),
			proving that \(x\in\argmin\varphi^{\nicefrac{\h}{\gamma}}\).
		\item\ref{thm:Moreau:lb}~
			It follows from \cref{thm:Moreau:eq} that if \(\varphi^{\nicefrac{\h}{\gamma}}\) is level bounded then so is \(\varphi\) on \(\interior\dom\h\).
			Conversely, suppose that there exists \(\alpha\in\R\) together with an unbounded sequence \(\seq{x^k}\subseteq[\varphi^{\nicefrac{\h}{\gamma}}\leq\alpha]\).
			Then, it follows from \cref{thm:Moreau:dom} that \(x^k\in\interior\dom\h\) for all \(k\), and in turn that for any \(k\) there exists \(\bar x^k\in\prox_{\gamma\varphi}^{\h}(x^k)\) which satisfies
			\(
				\varphi(\bar x^k)
			{}\leq{}
				\alpha-\D*(\bar x^k,x^k)
			\)
			by \cref{thm:Moreau:eq}.
			Local boundedness of \(\D\) with respect to the second variable (\cref{thm:Dhcoercive}) then ensures that \(\seq{\bar x^k}\subset\dom\h\) is not bounded, hence that \(\varphi^{\nicefrac{\h}{\gamma}}\) is not level bounded.
		\qedhere
		\end{proofitemize}
	\end{proof}
\end{prop}

It is apparent from \cref{thm:Moreau:dom} that the assertions of \cref{thm:Moreau:basic} cannot be extended outside of \(\interior\dom\h\).
This hindrance is also at the base of the requirement \(\range\prox_{\gamma\varphi}^{\h}\subseteq\interior\dom\h\); while guaranteed in case \(\varphi\) is convex (or if \(\dom\h=\R^n\)), it stands as a blanket assumption in Bregman optimization whenever these favorable conditions are not met.
It is also common in practice that the Bregman kernel \(\h\) is twice continuously differentiable, in which case we can easily derive subdifferential properties of the Bregman-Moreau envelope similarly to what done in \cite[\S3]{kan2012moreau}.

\begin{prop}[subdifferential properties of the Bregman-Moreau envelope]\label{thm:Moreau:subdifferential}%
	Let \(\h\) be a Legendre kernel with \(\h\in\C^2(\interior\dom\h)\) and \(\func{\varphi}{\R^n}{\Rinf}\).
	For every \(\gamma\in(0,\gamma_\varphi^{\h})\), the envelope \(\varphi^{\nicefrac{\h}{\gamma}}\) is strictly differentiable wherever it is differentiable.
	Moreover, for every \(x\in\interior\dom\h\)
	\begin{enumerate}
	\item
		\(
			\lip\varphi^{\nicefrac{\h}{\gamma}}(x)
		{}={}
			\max_{\bar x\in\prox_{\gamma\varphi}^{\h}(x)}\bigl\|
				\tfrac1\gamma\nabla^2\h(x)(x-\bar x)
			\bigr\|
		\);
	\item
		\(
			\partial\varphi^{\nicefrac{\h}{\gamma}}(x)
		{}={}
			\partial_B\varphi^{\nicefrac{\h}{\gamma}}(x)
		{}\subseteq{}
			\tfrac1\gamma\nabla^2\h(x)\bigl(x-\prox_{\gamma\varphi}^{\h}(x)\bigr)
		\).
	\end{enumerate}
	\begin{proof}
		We pattern the proof of \cite[Ex. 10.32]{rockafellar2011variational}.
		Let \(\mathcal U\subset\overline{\mathcal U}\subset\interior\dom\h\) be a bounded open set containing \(x\), and observe that by osc of \(\prox_{\gamma\varphi}^{\h}\) (\cref{thm:proxosc}) there exists a compact set \(\mathcal V\subseteq\dom\h\) such that
		\(
			-\varphi^{\h}(u)
		{}={}
			\max_{v\in\mathcal V}\Phi(u,v)
		\)
		for all \(u\in\mathcal U\), where \(\Phi(u,v)\coloneqq -\varphi(v)-\tfrac1\gamma\D_{\h}(v,u)\) is \(\C^1\) in \(u\), its derivatives depending continuously on \((u,v)\) with
		\(
			\nabla_u\Phi(u,v)
		{}={}
			\frac1\gamma\nabla^2\h(u)(v-u)
		\).
		In fact, the maxima are attained for \(v\in\prox_{\gamma\varphi}^{\h}(u)\).
		Function \(-\varphi^{\h}\) is thus \emph{lower-\(\C^1\)} in the sense of \cite[Def. 10.29]{rockafellar2011variational}, and all claims then follow from \cite[Thm. 10.31]{rockafellar2011variational}.
	\end{proof}
\end{prop}

		\subsection{Fixed points}
			In the Euclidean setting, it is well known that being a fixed point of the proximal mapping is an intermediate property between stationarity and global minimality, the three conditions being equivalent under convexity.
This is the rationale of the proximal point algorithm, which recursively applies the proximal mapping aiming at converging (possibly subsequentially and under boundedness assumptions) to one of its fixed points.
Major complications in the Bregman setting arise when the Legendre kernel has not full domain, starting from the fact that the boundary cannot contain fixed points, as is apparent from \cref{thm:Moreau:dom}.
Nevertheless, in the interior of \(\dom\h\) fixed points enjoy similar properties as in the Euclidean case and still play an important role in the algorithmic analysis.
As we will see in later subsections, single-valuedness of the proximal mapping at fixed points is of paramount importance for local regularity properties and thus deserves a dedicated definition.

\begin{defin}[nondegenerate fixed point]
	We say that \(x\in\R^n\) is a \DEF{fixed point} of the set-valued mapping \(\ffunc T{\R^n}{\R^n}\) if \(x\in T(x)\).
	We say that \(x\) is \DEF{nondegenerate} if \(T(x)=\set x\).
\end{defin}

Fortunately, degenerate fixed points are rarely encountered; in fact, with zero probability as validated in the following result.

\begin{lem}\label{thm:prox:nondegenerate}
	Let \(\func{\h}{\R^n}{\Rinf}\) be Legendre and \(\func{\varphi}{\R^n}{\Rinf}\) proper and lsc.
	If \(x\in\interior\dom\h\) is a fixed point of \(\prox_{\gamma\varphi}^{\h}\), then it is a nondegenerate fixed point of \(\prox_{\gamma\varphi}^{\h+\h'}\) for every strictly convex function \(\func{\h'}{\R^n}{\Rinf}\) differentiable on \(\interior\dom\h'\supseteq\interior\dom\h\).
	In particular, \(x\) is a nondegenerate fixed point of \(\prox_{\gamma'\varphi}^{\h}\) for every \(\gamma'\in(0,\gamma)\).
	\begin{proof}
		By definition, \(x\in\prox_{\gamma\varphi}^{\h}\) iff
		\(
			\varphi(z)
			{}+{}
			\tfrac{1}{\gamma}\D*(z,x)
		{}\geq{}
			\varphi(x)
		\)
		holds for any \(z\in\R^n\).
		Since \(\h+\h'\) is Legendre, \(\D*_{\h+\h'}=\D*+\D*_{\h'}\) on \(\dom\h\times\interior\dom\h\), and \(\D*_{\h'}(z,x)=0\) iff \(z=x\), apparently \(x\) is also a fixed point of \(\prox_{\gamma\varphi}^{\h+\h'}\), with
		\(
			\varphi(z)
			{}+{}
			\tfrac1\gamma\D*_{\h+\h'}(z,x)
		{}>{}
			\varphi(x)
		\)
		for all \(z\neq x\).
		Thus, if contrary to the claim there exists \(x'\in\prox_{\gamma\varphi}^{\h+\h'}(x)\setminus\set x\), then
		\[
			\varphi(x')
			{}+{}
			\tfrac1\gamma\D*(x',x)
		{}<{}
			\varphi(x)
		{}={}
			\varphi^{\nicefrac{(\h+\h')}{\gamma}}(x)
		{}={}
			\varphi(x')
			{}+{}
			\tfrac1\gamma\D*_{\h+\h'}(x',x)
		\]
		which is a contradiction.
	\end{proof}
\end{lem}

The following result shows that being a fixed point for the proximal map (for some stepsize \(\gamma\)) is actually a necessary condition for \emph{local} minimality in the interior of the domain of the Legendre kernel \(\h\).
An equivalence of local minimality for function \(\varphi\) and its envelope \(\varphi^{\nicefrac{\h}{\gamma}}\) at fixed points is also shown, in which nondegeneracy plays a key role.
We remark that this condition is necessary even in the Euclidean setting, as can be verified from the counterexample \(\varphi=\tfrac12\|{}\cdot{}\|^2+\indicator_{\set{0,1}}\) (see \cite[Fig. 3.1]{themelis2018proximal} and discussion therein).

\begin{thm}[equivalence of local minimality]\label{thm:Moreu:locmin}%
	The following hold for a Legendre kernel \(\func{\h}{\R^n}{\Rinf}\) and a proper, lsc, \(\h\)-prox-bounded function \(\func{\varphi}{\R^n}{\Rinf}\):
	\begin{enumerate}
	\item\label{thm:Moreau:locmin:fixed}
		if \(x_\star\in\interior\dom\h\) is a local minimum for \(\varphi\),
		then it is a nondegenerate fixed point of \(\prox_{\gamma\varphi}^{\h}\) for any \(\gamma\) small enough.
	\item\label{thm:Moreau:locmin:equiv}
		Conversely, any nondegenerate fixed point \(x_\star\) of \(\prox_{\gamma\varphi}^{\h}\) is a local minimum for \(\varphi\) iff it is a local minimum for \(\varphi^{\nicefrac{\h}{\gamma}}\).
		Moreover, equivalence of strong local minimality also holds provided that \(\h\) is strongly convex in a neighborhood of \(x_\star\).
	\end{enumerate}
	\begin{proof}
		\begin{proofitemize}
		\item\ref{thm:Moreau:locmin:fixed}
			For \(\gamma\in(0,\gamma_\varphi^{\h})\), let \(x_\gamma\in\prox_{\gamma\varphi}^{\h}(x_\star)\).
			It is easy to see that \(x_\gamma\to x_\star\) as \(\gamma\searrow0\) (cf. proof of \cite[Thm. 2.5]{kan2012moreau}).
			Local minimality of \(x_\star\) thus implies that \(\varphi(x_\gamma)\geq\varphi(x_\star)\) holds for \(\gamma\) small enough.
			Combined with \cref{thm:Moreau:eq} we obtain that
			\(
				\varphi(x_\gamma)+\tfrac1\gamma\D*(x_\gamma,x_\star)
			{}\leq{}
				\varphi(x_\star)
			{}\leq{}
				\varphi(x_\gamma)
			\)
			holds for \(\gamma\) small enough, hence that \(\D*(x_\gamma,x_\star)=0\) or, equivalently, \(x_\gamma=x_\star\).
		\item\ref{thm:Moreau:locmin:equiv}
			That (strong) local minimality for \(\varphi^{\nicefrac{\h}{\gamma}}\) implies that for \(\varphi\) follows from the fact that \(\varphi^{\nicefrac{\h}{\gamma}}\) ``supports'' \(\varphi\) at \(x_\star\), namely that \(\varphi^{\nicefrac{\h}{\gamma}}\leq\varphi\) and \(\varphi^{\nicefrac{\h}{\gamma}}(x_\star)=\varphi(x_\star)\) (\cref{thm:Moreau:eq}).
			Conversely, suppose that \(\h\) is \(\sigma_{\h,\mathcal U}\)-strongly convex in a neighborhood \(\mathcal U\) of \(x_\star\) for some \(\sigma_{\h,\mathcal U}\geq0\), and that there exists \(\mu\geq0\) such that
			\(
				\varphi(x)
			{}\geq{}
				\varphi(x_\star)
				{}+{}
				\tfrac\mu2\|x-x_\star\|^2
			\)
			for \(x\in\mathcal U\).
			Notice that in allowing \(\sigma_{\h,\mathcal U}=0\) and \(\mu=0\) we also cover nonstrong minimality and nonstrong convexity.
			Let \(\delta\coloneqq\tfrac12\min\set{\mu,\frac{\sigma_{\h,\mathcal U}}{2\gamma}}\geq0\), and note that \(\delta=0\) iff either \(\sigma_{\h,\mathcal U}\) or \(\mu\) is zero.
			To arrive to a contradiction, suppose that for all \(k\in\N_{\geq1}\) there exists \(x^k\in\ball{x_\star}{\nicefrac1k}\) such that
			\(
				\varphi^{\nicefrac{\h}{\gamma}}(x^k)
			{}<{}
				\varphi^{\nicefrac{\h}{\gamma}}(x_\star)
				{}+{}
				\tfrac\delta2\|x^k-x_\star\|^2
			{}={}
				\varphi(x_\star)
				{}+{}
				\tfrac\delta2\|x^k-x_\star\|^2
			\).
			Let \(\bar x^k\in\prox_{\gamma\varphi}^{\h}(x^k)\); since \(\prox_{\gamma\varphi}^{\h}\) is osc on \(\interior\dom\h\ni x_\star\) (cf. \cref{thm:facts}) and \(\prox_{\gamma\varphi}^{\h}(x_\star)=\set{x_\star}\), necessarily \(\bar x^k\to x_\star\) as \(k\to\infty\).
			We have
			\[
				\varphi(\bar x^k)
			{}\overrel*{\ref{thm:Moreau:eq}}{}
				\varphi^{\nicefrac{\h}{\gamma}}(x^k)
				{}-{}
				\tfrac1\gamma
				\D*(\bar x^k,x^k)
			{}\overrel*[\leq]{\ref{thm:hstrcvx}}{}
				\varphi^{\nicefrac{\h}{\gamma}}(x^k)
				{}-{}
				\tfrac{\sigma_{\h,\mathcal U}}{2\gamma}
				\|x^k-\bar x^k\|^2
			{}<{}
				\varphi(x_\star)
				{}+{}
				\tfrac\delta2\|x^k-x_\star\|^2
				{}-{}
				\tfrac{\sigma_{\h,\mathcal U}}{2\gamma}
				\|x^k-\bar x^k\|^2.
			\]
			By using the inequality
			\(
				\tfrac12\|a-c\|^2
			{}\leq{}
				\|a-b\|^2
				{}+{}
				\|b-c\|^2
			\)
			holding for any \(a,b,c\in\R^n\), we have
			\[
				\varphi(\bar x^k)
			{}<{}
				\varphi(x_\star)
				{}+{}
				\delta\|\bar x^k-x_\star\|^2
				{}+{}
				\bigl(\delta-\tfrac{\sigma_{\h,\mathcal U}}{2\gamma}\bigr)\|x^k-\bar x^k\|^2
			{}\leq{}
				\varphi(x_\star)
				{}+{}
				\tfrac\mu2\|\bar x^k-x_\star\|^2,
			\]
			where the last inequality follows from the definition of \(\delta\).
			Thus,
			\(
				\varphi(\bar x^k)
			{}<{}
				\varphi(x_\star)
				{}+{}
				\tfrac\mu2\|\bar x^k-x_\star\|^2
			\)
			for all \(k\in\N\), which contradicts (\(\mu\)-strong) local minimality of \(x_\star\) for \(\varphi\) (since \(\bar x^k\to x_\star\)).
			\qedhere
		\end{proofitemize}
	\end{proof}
\end{thm}

		\subsection{Local Euclidean reparametrization}
			The proximal mapping and the Moreau envelope in the Euclidean setting are clearly special instances of the more general Bregman variants, specifically corresponding to the choice \(\h=\tfrac12\|{}\cdot{}\|^2\) as Legendre kernel.
The converse, however, is not true even locally and in the convex case.
To see this, observe that while the Euclidean Moreau envelope preserves convexity, this is not at all the case for more general Bregman kernels even if they are of full domain, strongly convex and Lipschitz smooth.
Nevertheless, under a local strong convexity and Lipschitz differentiability assumption on \(\h\), yet with no requirement on its domain, it is possible to identify the Bregman proximal mapping and its Moreau envelope with a Euclidean forward-backward mapping and the corresponding envelope function, namely, the \DEF{forward-backward envelope} \cite{patrinos2013proximal,themelis2018forward}.
This result complements what was first observed in \cite{liu2017further}, namely that the Euclidean forward-backward envelope is, in fact, a Bregman-Moreau envelope.
The advantage of this identification will be revealed in the next subsections, where local differentiability results will be deduced with virtually no effort based on already established results in the Euclidean setting.
We remind the reader that given a decomposition \(\varphi=\FBsmooth+\FBnonsmooth\) with \(\FBsmooth\) continuously differentiable, the (Euclidean) forward-backward operator with stepsize \(\FBstepsize>0\) is
\begin{align}
\nonumber
	\FB x
{}={} &
	\argmin_{w\in\R^n}\set{
		\FBsmooth(x)
		{}+{}
		\innprod{\nabla\FBsmooth(x)}{z-x}
		{}+{}
		\FBnonsmooth(w)
		{}+{}
		\tfrac{1}{2\FBstepsize}\|w-x\|^2
	},
\shortintertext{%
	while the forward-backward envelope (FBE) is the associated value function, namely
}
\label{eq:FBE}
	\FBE*(x)
{}={} &
	\inf_{w\in\R^n}\set{
		\FBsmooth(x)
		{}+{}
		\innprod{\nabla\FBsmooth(x)}{w-x}
		{}+{}
		\FBnonsmooth(w)
		{}+{}
		\tfrac{1}{2\FBstepsize}\|w-x\|^2
	}.
\end{align}

\begin{thm}[local equivalence of Bregman-Moreau envelope and FBE]\label{thm:Moreau=FBE}%
	Let \(\func{\h}{\R^n}{\Rinf}\) be Legendre and \(\func{\varphi}{\R^n}{\Rinf}\) proper and lsc, and let \(\gamma\in(0,\gamma_\varphi^{\h})\) be fixed.
	Suppose further that \(\h\) is locally Lipschitz differentiable and locally strongly convex on \(\interior\dom\h\) (as is the case when \(\h\in\C^2\) with \(\nabla^2\h\succ0\) on \(\interior\dom\h\)) and that \(\range\prox_{\gamma\varphi}^{\h}\subseteq\interior\dom\h\).
	Then, for every compact set \(\mathcal U\subset\interior\dom\h\) there exist \(\bar\gamma>0\) and a convex compact set \(\mathcal V\) with \(\mathcal U\subseteq\mathcal V\subset\interior\dom\h\) such that for all \(\tilde\gamma\in(0,\bar\gamma)\) it holds that
	\begin{align}\label{eq:BMoreau=FBE}
		\varphi^{\nicefrac{\h}{\gamma}}=\FBE*
	\quad\text{and}\quad &
		\prox_{\gamma\varphi}^{\h}=\FB{}
		\text{ on \(\mathcal U\),}
	\shortintertext{where}
	\label{eq:Moreau:ftgt}
		\tilde f
		{}\coloneqq{}
			-\tfrac1\gamma\h+\tfrac{1}{2\tilde\gamma}\|{}\cdot{}\|^2
	\quad\text{and}\quad &
			\tilde g
		{}\coloneqq{}
			\varphi+\tfrac1\gamma\h-\tfrac{1}{2\tilde\gamma}\|{}\cdot{}\|^2+\indicator_{\mathcal V}.
	\end{align}
	Moreover, \(\tilde g\) is proper, lsc, and prox-bounded (in the Euclidean sense) with \(\gamma_{\tilde g}=\infty\), and \(\tilde f\) is \(L_{\tilde f}\)-Lipschitz-differentiable on \(\mathcal V\) with \(\tilde\gamma<\nicefrac{1}{L_{\tilde f}}\).
	\begin{proof}
		Let
		\(
			\mathcal V\coloneqq\conv(\mathcal U\cup\prox_{\gamma\varphi}^{\h}(\mathcal U))
		\)
		and observe that \(\mathcal V\) is a convex compact subset of \(\interior\dom\h\), as it follows from osc and local boundedness of \(\prox_{\gamma\varphi}^{\h}\) (\cref{thm:facts}) and the fact that compactness is preserved by the convex hull.
		It follows from the assumptions that there exist \(L_{\h,\mathcal V}\geq\sigma_{\h,\mathcal V}>0\) such that \(\h\) is \(L_{\h,\mathcal V}\)-smooth and \(\sigma_{\h,\mathcal V}\)-strongly convex on \(\mathcal V\).
		Define
		\(
			\bar\gamma
		{}\coloneqq{}
			\tfrac{2\gamma}{L_{\h,\mathcal V}}
		\),
		let \(\tilde\gamma\in(0,\bar\gamma)\) be fixed, and let \(\tilde f\) and \(\tilde g\) be as in \eqref{eq:Moreau:ftgt}.
		Clearly, \(\varphi=\tilde f+\tilde g\) on \(\mathcal V\).
		Moreover,
		\[
			\varphi(z)+\tfrac{1}{\gamma}\D*(z,x)+\indicator_{\mathcal V}(z)
		{}={}
			\tilde f(x)
			{}+{}
			\innprod{\nabla\tilde f(x)}{z-x}
			{}+{}
			\tilde g(z)
			{}+{}
			\tfrac{1}{2\tilde\gamma}\|z-x\|^2
		\quad
			\forall x\in\mathcal U,z\in\R^n.
		\]
		Since \(\prox_{\gamma\varphi}^{\h}(\mathcal U)\subseteq\mathcal V\),
		\[\ifams\mathtight[0.8]\fi
			\varphi^{\nicefrac{\h}{\gamma}}(x)
		{}={}
			\min_{z\in\mathcal V}\set{
				\varphi(z)+\tfrac{1}{\gamma}\D*(z,x)
			}
		{}={}
			\min_{z\in\R^n}\set{
				\tilde f(x)
				{}+{}
				\innprod{\nabla\tilde f(x)}{z-x}
				{}+{}
				\tilde g(z)
				{}+{}
				\tfrac{1}{2\tilde\gamma}\|z-x\|^2
			}
		{}={}
			\FBE*(x).
		\]
		Similarly, by considering the minimizers it is apparent that \(\ifams\mathtight[0.8]\fi\prox_{\gamma\varphi}^{\h}(x)=\FB x\) for any \(x\in\mathcal U\).
		Notice that \(\tilde g\) is proper, lsc and with bounded domain, hence its claimed prox-boundedness.
		Moreover,
		\[
			\left(
				\tfrac{1}{\tilde\gamma}-\tfrac{L_{\h,\mathcal V}}{\gamma}
			\right)\|x-y\|^2
		{}\leq{}
			\innprod{\nabla\tilde f(x)-\nabla\tilde f(y)}{x-y}
		{}\leq{}
			\left(
				\tfrac{1}{\tilde\gamma}-\tfrac{\sigma_{\h,\mathcal V}}{\gamma}
			\right)\|x-y\|^2
		\]
		for every \(x,y\in\mathcal V\).
		Therefore, \(\tilde f\) is \(L_{\tilde f}\)-smooth on \(\mathcal V\), with
		\(\ifams\mathtight[0.95]\fi
			L_{\tilde f}
		{}={}
			\max\set{
				\left|
					\tfrac{1}{\tilde\gamma}-\tfrac{L_{\h,\mathcal V}}{\gamma}
				\right|,~
				\left|
					\tfrac{1}{\tilde\gamma}-\tfrac{\sigma_{\h,\mathcal V}}{\gamma}
				\right|
			}
		\).
		Since
		\(
			\gamma
		{}<{}
			\bar\gamma
		{}={}
			\tfrac{2\gamma}{L_{\h,\mathcal V}}
		\),
		it follows that
		\(
			-\tfrac{1}{\tilde\gamma}
		{}<{}
			\tfrac{1}{\tilde\gamma}-\tfrac{L_{\h,\mathcal V}}{\gamma}
		{}\leq{}
			\tfrac{1}{\tilde\gamma}-\tfrac{\sigma_{\h,\mathcal V}}{\gamma}
		{}<{}
			\tfrac{1}{\tilde\gamma}
		\),
		proving that \(\tilde\gamma<\nicefrac{1}{L_{\tilde f}}\).
	\end{proof}
\end{thm}

		\subsection{First- and second-order properties}
			Although strict continuity ensures almost everywhere differentiability, with mild additional assumptions the Bregman-Moreau envelope can be shown to be continuously differentiable around critical points.
Thanks to the local equivalence shown in \Cref{thm:Moreau=FBE}, these requirements are the same as those ensuring similar properties in the Euclidean case.
These amount to prox-regularity, a condition which was first proposed in \cite{poliquin1996proxregular} and which has been recently extended to its \(\h\)-relative version in \cite{laude2020bregman}.

\begin{defin}[prox-regularity]\label{def:proxreg}%
	A function \(\func{\varphi}{\R^n}{\Rinf}\) is prox-regular at \(\bar x\) for \(\bar v\in\partial\varphi(\bar x)\) if it is locally lsc at \(\bar x\) and there exists \(r,\varepsilon>0\) such that
	\begin{equation}\label{eq:proxreg}
		\varphi(x')\geq\varphi(x)+\innprod{v}{x'-x}-\tfrac r2\|x'-x\|^2
	\end{equation}
	holds for all
	\(
		x,x'\in\ball{\bar x}{\varepsilon}
	\)
	and
	\(
		(x,v)\in\graph\partial\varphi
	\)
	with
	\(
		v\in\ball{\bar v}{\varepsilon}
	\)
	and
	\(
		\varphi(x)\leq\varphi(\bar x)+\varepsilon
	\).
\end{defin}

Differentiability property of the Bregman-Moreau envelope have been studied in \cite{teboulle1992entropic,bauschke2005joint} for jointly convex Bregman distances in the convex setting, and a similar analysis for the ``right'' envelope is provided in \cite{bauschke2018regularizing}.
For a nonconvex function \(\varphi\), \cite{kan2012moreau} shows global continuous differentiability of \(\varphi^{\nicefrac{\h}{\gamma}}\) under a global convexity assumption on \(\h+\gamma\varphi\).
What we provide next is instead a local result that requires local properties of \(\varphi\) around nondegenerate fixed points of the proximal map.
We remark that after our first submission of the paper a similar result appeared in \cite{laude2020bregman} in a more general setting.
We however offer our alternative proof as a means of emphasizing the favorable theoretical implications of the Euclidean equivalence stated in \cref{thm:Moreau=FBE}, which, ultimately, will lead us to the second-order result of \cref{thm:Moreau:2nd} which is instead novel.

\begin{thm}[continuous differentiability of the Bregman-Moreau envelope]\label{thm:Moreau:1st}%
	Suppose that the assumptions of \cref{thm:Moreau=FBE} hold, and let \(x_\star\) be a nondegenerate fixed point of \(\prox_{\gamma\varphi}^{\h}\).
	If \(\varphi\) is prox-regular at \(x_\star\) for \(v=0\), then there exists a neighborhood \(\mathcal U\) of \(x_\star\) on which the following statements are true:
	\begin{enumerate}
	\item\label{thm:proxLip}
		\(\prox_{\gamma\varphi}^{\h}\) is Lipschitz continuous (hence single valued);
	\item\label{thm:Moreau:C1}
		\(\varphi^{\nicefrac{\h}{\gamma}}\in\C^1\) with
		\(
			\nabla\varphi^{\nicefrac{\h}{\gamma}}(x)
		{}={}
			\tfrac1\gamma\nabla^2\h(x)\bigl(
				x-\prox_{\gamma\varphi}^{\h}(x)
			\bigr)
		\).
	\end{enumerate}
	\begin{proof}
			For any compact neighborhood \(\mathcal U\subset\interior\dom\h\) of \(x_\star\) we may invoke \cref{thm:Moreau=FBE} and identify \(\varphi^{\nicefrac{\h}{\gamma}}\) with the Euclidean FBE \(\FBE*\) on \(\mathcal U\) and \(\prox_{\gamma\varphi}^{\h}\) with
			\(\ifams\mathtight[0.95]\fi
				T\coloneqq\FB{}
			\),
			for some \(\tilde\gamma>0\) and \(\tilde f\) and \(\tilde g\) as in \eqref{eq:Moreau:ftgt}.
			It follows from \cite[Ex. 13.35]{rockafellar2011variational} and the continuous differentiability of \(\h\) that \(\tilde g\) is prox-regular at \(x_\star\) for \(-\nabla\tilde f(x_\star)\).
			Since \(\tilde f\) is \(\C^2\) around \(x_\star\) and \(T(x_\star)=\prox_{\gamma\varphi}^{\h}(x_\star)=\set{x_\star}\) by assumption, the setting of \cite[Thm. 4.7]{themelis2018forward} is satisfied\footnote{\label{footnote:Gamma}%
				The requirement \(\gamma\in(0,\Gamma(x_\star))\) in \cite[Thm. 4.7]{themelis2018forward} is only needed for ensuring that \(T(x_\star)=\set{x_\star}\) is a singleton (which here follows by assumption), and consequently that strict inequality in equation (4.4) therein holds for \(x\neq x_\star\).%
			}
			and thus \(T\) is Lipschitz continuous around \(x_\star\) and the Euclidean FBE \(\FBE*\) is \(\C^1\) around \(x_\star\) with
			\begin{equation}
				\nabla\varphi^{\nicefrac{\h}{\gamma}}(x)
			{}={}
				\nabla\FBE*(x)
			{}={}
				\tilde\gamma^{-1}
				\bigl[
					\I-\tilde\gamma\nabla^2\!\tilde f(x)
				\bigr]
				\bigl(
					x-T(x)
				\bigr)
			{}={}
				\tfrac1\gamma\nabla^2\h(x)
				\bigl(
					x-\prox_{\gamma\varphi}^{\h}(x)
				\bigr),
			\end{equation}
			which completes the proof.
	\end{proof}
\end{thm}

			We conclude the section with a second-order analysis \emph{at} (as opposed to around) fixed points.
In the same spirit of \cref{thm:Moreau:1st}, we will exploit the local identities assessed in \cref{thm:Legendre} to reduce the problem to the known Euclidean setting.
We will invoke the notion of epi-differentiability, the same requirement ensuring the existence of the Hessian for the Euclidean Moreau envelope in \cite{poliquin1995second}.
We refer the reader to \cite{poliquin1995second,rockafellar2011variational} for an extensive discussion on these generalized differentiability notions.

\begin{thm}[twice differentiability of the Bregman-Moreau envelope]\label{thm:Moreau:2nd}%
	Additionally to the assumptions of \cref{thm:Moreau:1st}, suppose that \(\varphi\) is (strictly) twice epi-differentiable at \(x_\star\) for \(v=0\), with generalized quadratic second-order epi-derivative.
	Then,
	\begin{enumerate}
	\item\label{thm:prox:C1}
		\(\prox_{\gamma\varphi}^{\h}\) is (strictly) differentiable at \(x_\star\);
	\item\label{thm:prox:C1symm}
		\(\prox_{\gamma\varphi}^{\h}\circ\nabla\conj{\h}\) is (strictly) differentiable at \(\nabla\h(x_\star)\) with symmetric and positive semidefinite Jacobian;
	\item\label{thm:Moreau:C2}
		\(\varphi^{\nicefrac{\h}{\gamma}}\) is (strictly) twice differentiable at \(x_\star\) with symmetric Hessian
		\[
			\nabla^2\varphi^{\nicefrac{\h}{\gamma}}(x_\star)
		{}={}
			\tfrac1\gamma\nabla^2\h(x_\star)[x_\star-J\prox_{\gamma\varphi}^{\h}(x_\star)].
		\]
	\end{enumerate}
	\begin{proof}
		As shown in the proof of \cref{thm:Moreau:1st}, for some \(\tilde\gamma>0\) and with \(\tilde f\) and \(\tilde g\) as in \eqref{eq:Moreau:ftgt} we may identify \(\varphi^{\nicefrac{\h}{\gamma}}\) with the Euclidean FBE \(\FBE*\) around \(x_\star\) and \(\prox_{\gamma\varphi}^{\h}\) with the Euclidean proximal gradient operator \(\FB{}\).
		It follows from \cite[Ex.s 13.18 and 13.25]{rockafellar2011variational} and the continuous differentiability of \(\h\) that \(\tilde g\) is prox-regular and (strictly) twice epi-differentiable at \(x_\star\) for \(-\nabla\tilde f(x_\star)\), with generalized quadratic second-order epi-derivative.
		With the same arguments as in the proof of \cref{thm:FBEC1} (cf. \cref{footnote:Gamma}), the setting of \cite[Thm. 4.10]{themelis2018forward} is satisfied and therefore:
		\begin{enumerate}[resume]
		\item
			\(\FB{}=\prox_{\gamma\varphi}^{\h}\) is (strictly) differentiable at \(x_\star\), which is assertion \ref{thm:prox:C1}.
		\item
			\(\FBE*=\varphi^{\nicefrac{\h}{\gamma}}\) is (strictly) twice differentiable at \(x_\star\) with symmetric Hessian
			\[
				\tfrac{1}{\tilde\gamma}\bigl[
					\I-\tilde\gamma\nabla^2\tilde f(x_\star)
				\bigr]\left(
					x_\star-J\left(\FB{x_\star}\right)
				\right)
			{}={}
				\tfrac1\gamma\nabla^2\h(x_\star)
				\bigl(
					x_\star-J\prox_{\gamma\varphi}^{\h}(x_\star)
				\bigr),
			\]
			which is assertion \ref{thm:Moreau:C2}.
		\item
			\(\prox_{\tilde\gamma\tilde g}\) is (strictly) differentiable at \(\Fw{x_\star}\) with symmetric Jacobian.
			We have
			{\mathtight[0.95]\begin{align*}
				\prox_{\tilde\gamma\tilde g}(s)
			{}={} &
				\argmin_{w\in\mathcal V}\set{
					\varphi(w)+\tfrac1\gamma\h(w)-\tfrac{1}{2\tilde\gamma}\|w\|^2
					{}+{}
					\tfrac{1}{2\tilde\gamma}\|w-s\|^2
				}
			\\
			{}={} &
				\argmin_{w\in\mathcal V}\set{
					\varphi(w)+\tfrac1\gamma\h(w)
					{}-{}
					\tfrac1\gamma\innprod{\tfrac{\gamma}{\tilde\gamma}s}{w}
				}
			{}={}
				\argmin_{w\in\mathcal V}\set{
					\varphi(w)
					{}+{}
					\tfrac1\gamma\D*\bigl(w,\nabla\conj{\h}(\tfrac{\gamma}{\tilde\gamma}s)\bigr)
				}.
			\end{align*}}%
			Since \(\prox_{\gamma\varphi}^{\h}(x)\) is contained in \(\mathcal V\) for all points close to \(x_\star\), as it is apparent from the definition of \(\mathcal V\) in the proof of \cref{thm:Moreau=FBE}, for all points \(s\) such that \(\tfrac{\gamma}{\tilde\gamma}s\) is close to \(\nabla\h(x_\star)\) the above formula coincides with \(\prox_{\gamma\varphi}\circ\nabla\conj{\h}(\tfrac{\gamma}{\tilde\gamma}s)\).
			The proof of assertion \ref{thm:prox:C1symm} now follows by observing that \(\Fw{x_\star}=\frac{\tilde\gamma}{\gamma}\nabla\h(x_\star)\).
		\qedhere
		\end{enumerate}
	\end{proof}
\end{thm}

	\section{Bregman forward-backward mapping and forward-backward envelope}\label{sec:FBS}
		As a last step towards the algorithm presented in the next section, we here analyze its main building block, the Bregman forward-backward operator.
We thus go back to the composite minimization setting of the investigated problem, here stated again for the reader's convenience
\[\tag{P}\label{eq:P}
	\minimize\varphi(x)\equiv f(x)+g(x)
\quad\stt{}
	x\in\overline C,
\]
and which will be addressed under the following assumptions.

\begin{ass}[requirements for composite minimization \eqref{eq:P}]\label{ass:basic}%
	The following hold:
	\begin{enumeratass}
	\item\label{ass:h}
		\(\func h{\R^n}{\Rinf\coloneqq\R\cup\set\infty}\) is a Legendre kernel with \(\interior\dom h=C\);
	\item\label{ass:f}
		\(\func f{\R^n}{\Rinf}\) is \(L_f\)-smooth relative to \(h\);
	\item\label{ass:g}
		\(\func g{\R^n}{\Rinf}\) is proper and lower semicontinuous (lsc);
	\item\label{ass:phi}
		\(\argmin\set{\varphi(x)}[x\in\overline C]\neq\emptyset\);
	\item\label{ass:rangeT}
		\(\range\T\subseteq C\) for \(\gamma\in(0,\nicefrac{1}{L_f})\), where
		\begin{equation}\label{eq:T}
			\T(x)
		{}\coloneqq{}
			\argmin_{z\in\R^n}\set{
				f(x)+\innprod{\nabla f(x)}{z-x}+g(z)+\tfrac1\gamma\D_h(z,x)
			}
		\end{equation}
		is the \DEF{Bregman forward-backward} mapping.
	\end{enumeratass}
\end{ass}

Our approach hinges on two analogies, one based on the local equivalence of the Bregman proximal map and the Euclidean forward-backward mapping given in \cref{thm:Moreau=FBE} and particularly useful for asymptotic analyses, and the other one based on the equivalence of forward-backward and proximal mappings in the Bregman setting.
The latter identity, which we show next in \cref{thm:Legendre}, leads to a simpler analysis of the Bregman forward-backward mapping in allowing us to disregard the decomposition \(f+g\) to solely focus on the cost function \(\varphi\).
Most importantly, it enables the possibility to make use of the Bregman-Moreau envelope in the algorithmic analysis, whence the thorough study carried out in the previous section will be heavily exploited.
In this perspective, in the same spirit of the Bregman-Moreau envelope and its relation with the Bregman proximal mapping, we construct an ``envelope'' for the Bregman forward-backward operator by considering the value function associated to the minimization problem \eqref{eq:T} defining \(\T\).
We define the \DEF{Bregman forward-backward envelope (BFBE)} as the function \(\func{\FBE}{\R^n}{\Rinf}\) given by
\begin{equation}\label{eq:BFBE}
	\FBE(x)
{}\coloneqq{}
	\inf_{z\in\R^n}\set{
		f(x)+\innprod{\nabla f(x)}{z-x}
		{}+{}
		g(z)
		{}+{}
		\tfrac1\gamma\D(z,x)
	}.
\end{equation}
As anticipated, in the proximal point reinterpretation of \(\T\) the BFBE corresponds to the Bregman-Moreau envelope, as shown next.

\begin{thm}[equivalence of forward-backward and proximal point mappings]\label{thm:Legendre}%
	Suppose that \(f\) is \(L_f\)-relatively smooth with respect to a Legendre kernel \(h\).
	Then, for every \(\gamma\in(0,\nicefrac{1}{L_f})\) the function \(\hat h\coloneqq\frac h\gamma-f\) (with the convention \(\infty-\infty=\infty\)) is a Legendre kernel.
	Moreover,
	\begin{align}
	\label{eq:M=PPM}
		\varphi(z)+\D*(z,x)
	{}={} &
		f(x)+\innprod{\nabla f(x)}{z-x}
		{}+{}
		g(z)
		{}+{}
		\tfrac1\gamma\D(z,x)
	\shortintertext{%
		holds for any \((z,x)\in\R^n\times\interior\dom h\), and in particular
	}
	\label{eq:Tprox}
		\T(x)
	{}={} &
		\prox_{\varphi}^{\nicefrac h\gamma-f}(x)
	{}={}
		\prox_{\gamma\varphi}^{h-\gamma f}(x)
	\shortintertext{%
		and
	}
	\label{eq:BFBE=Moreau}
		\FBE(x)
	{}={} &
		\varphi^{\nicefrac h\gamma-f}(x).
	\end{align}
	\begin{proof}
		Let \(C=\interior\dom h\).
		Observe that
		\(
			\hat h
		{}={}
			\tfrac{1-\gamma L_f}{\gamma}h
			{}+{}
			L_fh-f
		\)
		on \(\dom\hat h=\dom h\); since \(L_fh-f\) is convex on \(C\) by assumption, \(1\)-coercivity and (essential) strict convexity on \(C\) of \(\hat h\) follow from the similar properties of \(h\).
		We now show essential smoothness; clearly, \(\hat h\) is differentiable on \(C\) with \(\nabla\hat h=\tfrac1\gamma\nabla h-\nabla f\).
		To arrive to a contradiction, suppose that there exists a sequence \(\seq{x^k}\subset C\) converging to a boundary point \(x_\star\) of \(C\) and such that \(\sup_{k\in\N}\|\nabla\hat h(x^k)\|<\infty\).
		By possibly extracting a subsequence, we may assume that \(\nabla h(x^k)/\|\nabla h(x^k)\|\to v\) for some unitary vector \(v\).
		For every \(y\in C\), since
		\(
			\nabla(L_fh-f)(x^k)
		{}={}
			\nabla\hat h(x^k)
			{}-{}
			\tfrac{1-\gamma L_f}{\gamma}\nabla h(x^k)
		\),
		it holds that
		\begin{equation}\label{eq:neginnprod}
			\innprod{
				\nabla(L_fh-f)(x^k)
				{}-{}
				\nabla(L_fh-f)(y)
			}{
				x^k
				{}-{}
				y
			}
		{}\leq{}
			c_y
			{}-{}
			\tfrac{1-\gamma L_f}{\gamma}
			\innprod{
				\nabla h(x^k)
				{}-{}
				\nabla h(y)
			}{
				x^k
				{}-{}
				y
			},
		\end{equation}
		where
		\(
			c_y
		{}\coloneqq{}
			\sup_{k\in\N}\innprod{\nabla\hat h(x^k)-\nabla\hat h(y)}{x^k-y}
		\)
		is a finite quantity.
		Moreover, since \(\|\nabla h(x^k)\|\to\infty\),
		\[
			0
		{}\leq{}
			\tfrac{1}{\|\nabla h(x^k)\|}
			\innprod{
				\nabla h(x^k)
				{}-{}
				\nabla h(y)
			}{
				x^k
				{}-{}
				y
			}
		{}\to{}
			\innprod{v}{x_\star-y}
		\quad
			\text{as }
			k\to\infty,
		\]
		and from the arbitrarity of \(y\in C\) we conclude that
		\(
			v
		{}\in{}
			\set{u}[
				\innprod{u}{x_\star-y}\geq0
				~\forall y\in C
			]
		\).
		Since \(C\) is open, \(\ball{x_\star}{\varepsilon}\cap C\neq\emptyset\) for any \(\varepsilon>0\), and in particular there exists \(y\in C\) such that \(\innprod{v}{x_\star-y}\gneqq0\).
		Pluggin this \(y\) in \eqref{eq:neginnprod} yields
		\[
			\innprod{
				\nabla(L_fh-f)(x^k)
				{}-{}
				\nabla(L_fh-f)(y)
			}{
				x^k
				{}-{}
				y
			}
		{}\leq{}
			c_y
			{}-{}
			\tfrac{1-\gamma L_f}{\gamma}
			\|\nabla h(x^k)\|
			\innprod{
				\tfrac{
					\nabla h(x^k)
					{}-{}
					\nabla h(y)
				}{
					\|\nabla h(x^k)\|
				}
			}{
				x^k
				{}-{}
				y
			}
		{}\to{}
			-\infty,
		\]
		contradicting convexity of \(L_fh-f\) on \(C\).
		Therefore, \(\tfrac1\gamma h-f\) is a Legendre kernel, and in particular the right-hand side in \eqref{eq:M=PPM} is well defined, with equality therein and consequent validity of \eqref{eq:Tprox} and \eqref{eq:BFBE=Moreau} of immediate verification.
	\end{proof}
\end{thm}

In light of the equivalence of \cref{thm:Legendre}, properties of the Bregman proximal mapping and Bregman-Moreau envelope can directly be imported.

\begin{prop}\label{thm:T}
	Suppose that \cref{ass:basic} holds and let \(\gamma\in(0,\nicefrac{1}{L_f})\).
	Then,
	\begin{enumerate}
	\item\label{thm:Tprox}%
		\(
			\T(x)
		{}={}
			\prox_{\varphi}^{\nicefrac h\gamma - f}(x)
		{}={}
			\prox_{\gamma g}^h\bigl(\nabla \conj h(\nabla h(x)-\gamma\nabla f(x))\bigr)
		\)
		for every \(x\in C\);
	\item\label{thm:T:dom}%
		\(\dom\FBE=\dom\T=C\);
	\item
		\(\range\T\subseteq\dom\varphi\cap C\);
	\item\label{thm:Tosc}%
		\(\T\) is locally bounded, compact-valued, and osc on \(C\);
	\item\label{thm:FBEC0}%
		\(\FBE\) is continuous on \(C\); in fact, locally Lipschitz if so are \(\nabla f\) and \(\nabla h\);
	\item\label{thm:nondegenerate}
		if \(x\in\T(x)\), then \(\renewcommand\Tstepsize{\gamma'}\T(x)=\set x\) for every \(\gamma'\in(0,\gamma)\).
	\end{enumerate}
	\begin{proof}
		It suffices to show the second equality in assertion \ref{thm:Tprox}, while all other claims follow from the similar properties of the Bregman proximal mapping and Bregman-Moreau envelope in light of the equivalence of \cref{thm:Legendre}.
		By expanding the Bregman distance and discarding constant terms in \eqref{eq:T}, one has
		\[
			\T(x)
		{}={}
			\argmin_z\set{
				g(z)
				{}+{}
				\tfrac1\gamma\bigl[
					h(z)-\innprod{\nabla h(x)-\gamma\nabla f(x)}{z-x}
				\bigr]
			}
		{}={}
			\argmin_z\set{
				g(z)
				{}+{}
				\tfrac1\gamma\D(z,\bar z)
			}
		\]
		for
		\(
			\bar z=\nabla\conj h(\nabla h(x)-\gamma\nabla f(x))
		\),
		owing to the identity
		\(
			\nabla h(x)-\gamma\nabla f(x)
		{}={}
			\nabla h(\bar z)
		\)
		(cf. \cref{thm:inv1}), hence the claim.
	\end{proof}
\end{prop}

		The next two results characterize the fundamental relationship between the Bregman forward-backward envelope \(\FBE\) and the original function \(\varphi\) that are essential to analyze the convergence of the Bregman forward-backward scheme that will be given in \Cref{sec:Algorithm}.

\begin{prop}[relation between \(\varphi\) and \(\FBE\)]%
	Suppose that \Cref{ass:basic} holds and let \(\gamma\in(0,\nicefrac{1}{L_f})\) be fixed.
	Then,
	\begin{enumerate}
	\item\label{thm:leq}
		\(\FBE(x)\leq\varphi(x)\) for all \(x\in C\), with equality holding iff \(x\in\T(x)\);
	\item\label{thm:geq}
		\(
			\tfrac{1-\gamma L_f}{\gamma}\D_h(\bar x,x)
		{}\leq{}
			\FBE(x)-\varphi(\bar x)
		{}\leq{}
			\tfrac{1+\gamma L_f}{\gamma}\D_h(\bar x,x)
		\)
		for all \(x\in C\) and \(\bar x\in \T(x)\);
	\item\label{thm:inf}
		\(\inf\FBE=\inf_C\varphi\) and \(\argmin\FBE=\argmin_C\varphi\);
	\item\label{thm:LB}
		\(\FBE\) is level bounded iff \(\varphi\) is level bounded on \(C\).
	\end{enumerate}
	\begin{proof}
		All the claims follow from \cref{thm:Legendre,thm:Moreau:basic} together with the fact that
		\(
			\frac{1-\gamma L_f}{\gamma}\D
		{}\leq{}
			\D_{\hat h}
		{}\leq{}
			\frac{1+\gamma L_f}{\gamma}\D
		\)
		for \(\hat h=\frac1\gamma h-f\) (owing to convexity of \(L_fh\pm f\)).
	\end{proof}
\end{prop}

\begin{thm}[equivalence of local minimality]\label{thm:locmin}%
	Suppose that \Cref{ass:basic} holds.
	Then,
	\begin{enumerate}
	\item\label{thm:locmin:fixed}
		if \(x_\star\in C\) is a local minimum for \(\varphi\), then it is a nondegenerate fixed point of the proximal-gradient mapping \(\T\) for any \(\gamma\) small enough.
	\item\label{thm:locmin:equiv}
		Conversely, any nondegenerate fixed point \(x_\star\) of \(\T\) is a local minimum for \(\varphi\) iff it is a local minimum for \(\varphi^{\nicefrac h\gamma}\).
		Moreover, equivalence of strong local minimality also holds provided that \(h-\gamma f\) is strongly convex in a neighborhood of \(x_\star\).
	\end{enumerate}
	\begin{proof}
		Assertion \ref{thm:locmin:equiv} directly follows from \cref{thm:Moreau:locmin:equiv} in light of the equivalence of \cref{thm:Legendre}.
		Suppose now that \(x_\star\in C\) is a local minimum for \(\varphi\).
		Then, there exists \(\bar\gamma>0\) such that \(\prox_\varphi^{\nicefrac{h}{\gamma'}}(x_\star)=\set{x_\star}\) for all \(\gamma'\in(0,\bar\gamma)\), as shown in \cref{thm:Moreau:locmin:fixed}.
		We now claim that \(x_\star\) is a nondegenerate fixed point of \(\T\) for \(\gamma\in(0,\frac{\bar\gamma}{1+\bar\gamma L_f})\).
		To see this, recall that \(\T=\prox_{\varphi}^{\frac1\gamma h-f}\) and observe that
		\(
			\tfrac1\gamma h-f
		{}={}
			\frac{1}{\gamma'}h
			{}+{}
			(L_fh-f)
		\)
		for \(\gamma'=\frac{\gamma}{1-\gamma L_f}\).
		By invoking \cref{thm:prox:nondegenerate} we conclude that \(x_\star\) is a nondegenerate fixed point of \(\T\) whenever \(\gamma\) is such that \(\gamma'\in(0,\bar\gamma)\), that is, for \(\gamma\in(0,\frac{\bar\gamma}{1+\bar\gamma L_f})\) as claimed. 
	\end{proof}
\end{thm}

\Cref{thm:Legendre} also implies through \cref{thm:Moreau=FBE} that the Bregman FBE \eqref{eq:FBE} is locally equivalent to its Euclidean version \eqref{eq:BFBE}.

\begin{thm}[local equivalence of Bregman and Euclidean FBE]\label{thm:EFBE}%
	Suppose that \Cref{ass:basic} holds and let \(\gamma<\nicefrac{1}{L_f}\) be fixed.
	Suppose further that \(h\) is locally Lipschitz differentiable and locally strongly convex on \(C\) (as is the case when \(h\in\C^2\) with \(\nabla^2h\succ0\) on \(C\)).
	Then, for every compact set \(\mathcal U\subset C\) there exist \(\bar\gamma>0\) and a compact convex set \(\mathcal V\subset C\) such that for all \(\tilde\gamma\in(0,\bar\gamma)\) it holds that
	\begin{align}\label{eq:BFBE=FBE}
		\FBE=\FBE*
	\quad\text{and}\quad &
		\T=\T*
		\text{ on \(\mathcal U\),}
	\shortintertext{where}
	\label{eq:ftgt}
		\tilde f
		{}\coloneqq{}
			f-\tfrac1\gamma h+\tfrac{1}{2\tilde\gamma}\|{}\cdot{}\|^2
	\quad\text{and}\quad &
			\tilde g
		{}\coloneqq{}
			g+\tfrac1\gamma h-\tfrac{1}{2\tilde\gamma}\|{}\cdot{}\|^2+\indicator_{\mathcal V}.
	\end{align}
	Moreover, \(\tilde g\) is proper, lsc, and prox-bounded (in the Euclidean sense) with \(\gamma_{\tilde g}=\infty\), and \(\tilde f\) is \(L_{\tilde f}\)-Lipschitz-differentiable on \(\mathcal V\) with \(\tilde\gamma<\nicefrac{1}{L_{\tilde f}}\).
\end{thm}

We remark that, differently from the case of \cref{thm:Moreau=FBE}, it is strong convexity of \(h-\gamma f\) that is required, while that of \(h\) is sufficient but not necessary for the purpose.
In fact, local strong convexity of \(-f\) would waive the need for a similar requirement on \(h\).

		\subsection{First- and second-order properties}
			This subsection
\ifams
	provides
\else
	lists
\fi
some (sub)differential properties of the BFBE and the Bregman forward-backward operator.
All the results fall as a direct consequence of the similar ones derived in the previous section.
We remind the reader that this passage hinges on the key equivalence assessed in \cref{thm:Legendre}, namely \(\T=\prox_{\varphi}^{\h}\) and \(\FBE=\varphi^{\h}\) for \(\h=\tfrac1\gamma h-f\).
In particular, the following result is a direct consequence of \cref{thm:Moreau:subdifferential}.
For the sake of a lighter notation, it is convenient to introduce the matrix-valued mapping (defined wherever it makes sense)
\begin{equation}\label{eq:Q}
	\Q(x)
{}\coloneqq{}
	\tfrac1\gamma\nabla^2h(x)-\nabla^2f(x).
\end{equation}

\begin{prop}[subdifferential properties of the BFBE]%
	Additionally to \Cref{ass:basic}, suppose that \(f,h\in\C^2(C)\).
	For every \(\gamma\in(0,\nicefrac{1}{L_f})\) the BFBE \(\FBE\) is strictly differentiable wherever it is differentiable.
	Moreover, for every \(x\in C\) and with \(\Q\) as in \eqref{eq:Q}
	\begin{enumerate}
	\item
		\(
			\lip\FBE(x)
		{}={}
			\max_{\bar x\in\T(x)}\bigl\|
				\Q(x)(x-\bar x)
			\bigr\|
		\);
	\item\label{thm:subdifferential}
		\(
			\partial\FBE(x)
		{}={}
			\partial_B\FBE(x)
		{}\subseteq{}
			\Q(x)\bigl(x-\T(x)\bigr)
		\).
	\end{enumerate}
\end{prop}

\begin{thm}[continuous differentiability of the BFBE]\label{thm:1st}%
	Suppose that \Cref{ass:basic} holds and that \(f,h\in\C^2\) with \(\nabla^2h\succ0\) on \(C\).
	Suppose further that \(g\) is prox-regular at a nondegenerate fixed point \(x_\star\) of \(\T\) for \(-\nabla f(x_\star)\).
	Then, there exists a neighborhood \(\mathcal U\) of \(x_\star\) on which the following statements are true:
	\begin{enumerate}
	\item\label{thm:TLip}
		\(\T\) is Lipschitz continuous (hence single valued);
	\item\label{thm:FBEC1}
		\(\FBE\in\C^1(\mathcal U)\) with \(\nabla\FBE(x)=\Q(x)(x-\T(x))\), where \(\Q\) is as in \eqref{eq:Q}.
	\end{enumerate}
	\begin{proof}
		We may invoke \cite[Ex. 13.35]{rockafellar2011variational} to infer that \(\varphi\) is prox-regular at \(x_\star\) for \(v=0\).
		In particular, the assumptions of \cref{thm:Moreau:1st} are satisfied for the kernel \(\h=\frac1\gamma h-f=(L_fh-f)+(\frac1\gamma-L_f)h\), and the proof then follows from \cref{thm:Moreau:1st} in light of \cref{thm:Legendre}.
	\end{proof}
\end{thm}

			Twice differentiability of the BFBE will play a key role in the asymptotic analysis of the \cref{alg:Bella} algorithm discussed in the next section, when directions of quasi-Newton type are considered (cf. \cref{thm:DM}).
The following result offers sufficient conditions ensuring this property at a fixed point of the Bregman forward-backward mapping \(\T\).

\begin{thm}[twice differentiability of \(\FBE\)]\label{thm:2nd}%
	Additionally to \Cref{ass:basic}, suppose that
	\begin{enumeratass}
	\item
		\(f\in\C^2(C)\) and \(\nabla^2f\) is [strictly] continuous around nondegenerate fixed point \(x_\star\) of \(\T\);
	\item
		\(h\in\C^2(C)\) with \(\nabla^2h\succ0\);
	\item
		\(g\) is prox-regular and [strictly] twice epi-differentiable at \(x_\star\) for \(-\nabla f(x_\star)\), with its second-order epi-derivative being generalized quadratic.
	\end{enumeratass}
	Then, with \(\Q\) as in \eqref{eq:Q}, \(\T\) is [strictly] differentiable at \(x_\star\), and \(\FBE\) is [strictly] twice differentiable at \(x_\star\) with symmetric Hessian
	\[
		\nabla^2\FBE(x_\star)
	{}={}
		\Q(x_\star)\bigl[\I-J\T(x_\star)],
	\]
	\begin{proof}
		It follows from \cite[Ex.s 13.18 and 13.25]{rockafellar2011variational} that \(\varphi\) is prox-regular and [strictly] twice epi-differentiable at \(x_\star\) for \(-\nabla\tilde f(x_\star)\), with generalized quadratic second-order epi-derivative.
		The assumptions of \cref{thm:Moreau:2nd} are thus satisfied for the kernel \(\h=\frac1\gamma h-f=(L_fh-f)+(\frac1\gamma-L_f)h\), and the proof then follows from \cref{thm:Moreau:2nd} in light of \cref{thm:Legendre}.
	\end{proof}
\end{thm}

	\section{The \algname\ algorithm}\label{sec:Algorithm}
		Having completed an in-depth analysis of the Bregman forward backward mapping, we are now ready to introduce a new algorithm based on this building block.
The purpose of the algorithm is to globalize the convergence of fast local method for solving problem \eqref{eq:P}, exclusively by means of calls to the forward-backward operator \(\T\).
Our motivation stems from the fact that in the nonsmooth (and nonconvex) setting of problem \eqref{eq:P} typical linesearch strategies such as Armijo- or Wolfe-type are not applicable, and convergence of fast methods such as Newton-type can thus be guaranteed only if starting close to a solution, a requirement that makes any such approach only conceptual and of no practical utility.
The \emph{Bregman envelope linesearch algorithm} (\cref{alg:Bella}) comes in response to this issue.

\begin{algorithm}[ht]
	\algcaption{\algname\ ({\tt B}regman {\sc\texttt e}nve{\sc\texttt l}ope {\sc\texttt l}inesearch {\sc\texttt a}lgorithm)}%
	\label{alg:Bella}
	\begin{algorithmic}[1]
\setlength\baselineskip{1.2\baselineskip}%
\Require
	\begin{tabular}[t]{@{}l@{}}
		with \cref{ass:basic} holding, select
		stepsize \(\gamma\in(0,\nicefrac{1}{L_f})\),
		initial point \(x^0\in C\),%
	\\
		\(\sigma\in(0,\frac{1-\gamma L_f}{\gamma})\),~
		tolerance \(\varepsilon>0\),~
		max number of backtrackings \(i_{\rm max}\in\N\cup\set\infty\)
	\end{tabular}
\Initialize
	\(k=0\)%
\State\label{state:BFBS}%
	choose \(\bar x^k\in\T(x^k)\)
\State\label{state:Stop}%
	{\bf if}~ \(\D(\bar x^k,x^k)\leq\varepsilon\) ~{\bf then}~ \Return \(\hat x\coloneqq\bar x^k\) ~{\bf end if}%
\State\label{state:d}%
	choose a direction \(d^k\in\R^n\) and set \(\tau_k=1\) and \(i_k=0\)%
\State\label{state:x+}%
	\(
		x^{k+1}
	{}={}
		(1-\tau_k)\bar x^k + \tau_k(x^k+d^k)
	\)
\If{~ \(\FBE(x^{k+1})\leq\FBE(x^k)-\sigma\D_h(\bar x^k,x^k)\) ~}\label{state:LS}
	\Comment{Linesearch passed}
	\State\label{state:nm}%
		\(k\gets k+1\)~ and go to \cref{state:BFBS}%
\ElsIf{~ \(i_k=i_{\rm max}\) ~}
	\Comment{Max \#backtrackings: do plain BFBS step}
		\State
			\(
				x^{k+1}
			{}={}
				\bar x^k
			\),~
			\(k\gets k+1\)~ and go to \cref{state:BFBS}%
\Else\label{state:backtrack}%
	\Comment{Linesearch failed: backtrack and retry}
	\State
		\(\tau_k\gets\nicefrac{\tau_k}{2}\),~
		\(i_k\gets i_k+1\)~ and go to \cref{state:x+}
\EndIf
\end{algorithmic}

\end{algorithm}%

Having assessed, under \cref{ass:basic}, the continuity property of the BFBE and the inequality \(\FBE(\bar x)\leq\FBE(x)-\tfrac{1-\gamma L_f}{\gamma}\D(\bar x,x)\) holding and for any \(\bar x\in\T(x)\), the algorithmic rationale is quite self-explanatory.
At each iteration, once a direction \(d\) has been selected according to an arbitrary user-defined criterion (ideally a ``fast'' direction, but virtually any choice works regardless), the candidate update direction \(x+d\) is ``pushed'' towards the forward-backward step \(\bar x\) until a descent inequality is satisfied on the BFBE.
This well definedness aspect will be better detailed in the dedicated \cref{sec:Finite}, where a qualitative measure of stationarity of the output point \(\hat x\) is also given.
The remaineder of the section is then devoted to asymptotic analyses, namely subsequential (\cref{thm:subseq}), global and linear (\cref{thm:global,thm:linear}), and ultimately superlinear convergence even to nonisolated local minima (\cref{thm:tau1}).
Lastly, a Dennis-Moré-type criterion in the flavor of classical quasi-Newton analysis is also linked to superlinear convergence to strong local minima (\cref{thm:DM}).

The oracle complexity of one iteration of \Cref{alg:Bella} is dictated by three operations: a Bregman forward-backward call at \cref{state:BFBS}, the computation of a direction \(d^k\) at \cref{state:d}, and the evaluation of the BFBE at \cref{state:LS}.
As discussed in \cref{sec:superlinear}, for instance, directions of quasi-Newton type only involve direct linear algebra operations on already available quantities.
Moreover, the evaluation of \(\FBE(x^{k+1})\) only requires one call to \(\T(x^{k+1})\) (which can be stored and reused at \cref{state:BFBS}), and consequently, apart from the freedom in choosing suitably inexpensive directions, each iteration requires one call to the forward-backward operator per backtracking trial at \cref{state:LS}.
A bound on the number of these calls can be imposed by selecting a finite threshold \(i_{\rm max}\), which however makes no difference from a theoretical standpoint.
We also remark that \cref{alg:Bella} includes known methods as special cases; by setting \(d^k=\bar x^k-x^k\) it reduces to the Bregman forward-backward algorithm given in \cite{bolte2018first} (the linesearch condition \eqref{eq:LS} is satisfied regardless of the stepsize \(\tau_k\) owing to \Cref{thm:leq,thm:geq}), while for the Euclidean kernel \(h=\tfrac12\|{}\cdot{}\|^2\) one obtains the {\tt PANOC} algorithm given in \cite{stella2017simple}.
It is also worth remarking that by considering \(f=0\) and suitably choosing \(h\), \cref{alg:Bella} offers a linesearch extension to any algorithm that can be interpreted as a Bregman proximal point scheme.

		\subsection{Well definedness and finite termination}\label{sec:Finite}
			As discussed in the beginning of the section, the rationale of the linesearch involved in \cref{alg:Bella} is a simple consequence of basic properties of the BFBE.
The next result validates this claim.

\begin{lem}[well definedness of {\cref{alg:Bella}}]\label{lem:wellDefin}%
	Let \Cref{ass:basic} hold, and let \(\gamma\in(0,\nicefrac{1}{L_f})\) and \(\sigma\in(0,\frac{1-\gamma L_f}{\gamma})\) be fixed.
	Then, for any \(x\in C\), \(\bar x\in\T(x)\setminus\set x\) and \(d\in\R^n\) there exists \(\bar\tau\in(0,1]\) such that for any \(\tau\in[0,\bar\tau]\) the point
	\(
		x_\tau^+
	{}\coloneqq{}
		(1-\tau)\bar x+\tau(x+d)
	\)
	satisfies
	\[
		\FBE(x_\tau^+)
	{}\leq{}
		\FBE(x)
		{}-{}
		\sigma\D(\bar x,x).
	\]
	In particular, since \(x_\tau^+\in C\), the iterates of \Cref{alg:Bella} are well defined with linesearch at \cref{state:LS} terminating after a finite number of backtrackings regardless of the choice of \(\seq{d^k}\) and whether or not a finite maximum number of backtrackings \(i_{\rm max}\) is imposed.
	\begin{proof}
		It follows from \cref{thm:leq,thm:geq} that the strict inequality
		\[
			\FBE(x')
		{}<{}
			\FBE(x)
			{}-{}
			\sigma\D(\bar x,x)
		\]
		holds for \(x'=\bar x\in C\).
		Continuity of the envelope \(\FBE\) and openness of its domain as asserted in \cref{thm:FBEC0} then ensure that there exists a neighborhood \(\mathcal U\) of \(\bar x\) such that the inequality remains valid for any \(x'\in\mathcal U\).
		The proof now follows by observing that \(x_\tau^+\to\bar x\) as \(\tau\searrow0\) for any \(d\in\R^n\), i.e., \(x_\tau^+\in\mathcal U\) for small enough \(\tau\).
		Therefore, the claimed inequality will be satisfied after a finite number of backtrackings, which is our desired result.
	\end{proof}
\end{lem}

Next, we offer an explicit bound on the number of iterations needed to satisfy the termination criterion at \cref{state:Stop}, and provide a qualitative measure of stationarity of the output vector \(\hat x\).
We remind that the number of calls to the forward-backward operator \(\T\) (corresponding to the number of backtrackings at \cref{state:backtrack}) can artificially be bounded by means of selecting a finite parameter \(i_{\rm max}\) at initialization.
By doing so, one obtains that \cref{alg:Bella} terminates with at most \(i_{\rm max}\) times the number of iterations many calls to \(\T\) (although this is likely a massively loose estimate when sensible directions are employed).

\begin{thm}[{iteration complexity of \Cref{alg:Bella}}]%
	Suppose that \Cref{ass:basic} holds.
	Then,
	\begin{enumerate}
	\item\label{thm:finite}
		\Cref{alg:Bella} terminates within
		\(
			k
		{}\leq{}
			\frac{\varphi(x^0)-\inf\varphi}{\sigma\varepsilon}
		\)
		iterations;
	\item\label{thm:suboptimal}
		if \(\dom h=\R^n\) and \(h\) is \(\sigma_{h,\mathcal U}\)-strongly convex and \(L_{h,\mathcal U}\)-Lipschitz differentiable on an open convex set \(\mathcal U\) that contains all the iterates \(x^k\) and \(\bar x^k\) (this being true if \(h\in\C^2\) with \(\nabla^2 h\succ0\)), then the point \(\hat x\) returned by the algorithm satisfies
		\(
			\dist(0,\hat\partial\varphi(\hat x))
		{}\leq{}
			\tfrac{1+\gamma L_f}{\gamma}
			\sqrt{
				\tfrac{2L_{h,\mathcal U}^2}{\sigma_{h,\mathcal U}}
				\varepsilon
			}
		\).
	\end{enumerate}
	\begin{proof}
		We begin by observing that for every \(k\in\N\) it holds that
		\begin{equation}\label{eq:LS}
			\FBE(\bar x^{k+1})
		{}\leq{}
			\FBE(x^{k+1})
		{}\leq{}
			\FBE(x^k)-\sigma\D_h(\bar x^k,x^k).
		\end{equation}
		The first inequality owes to \cref{thm:leq,thm:geq}, whereas the second one is apparent when the condition at \cref{state:LS} is satisfied, and follows from \cref{thm:leq,thm:geq} together with the fact that \(\sigma<\frac{1-\gamma L_f}{\gamma}\) otherwise (that is, when the maximum number of backtrackings is reached and the nominal step \(x^{k+1}=\bar x^k\) is taken).
		\begin{proofitemize}
		\item\ref{thm:finite}~
			By telescoping the inequality \eqref{eq:LS} over the first \(K>0\) iterations we have
			\begin{equation}\label{eq:telescope}
				\sigma
				\sum_{k=0}^{K-1}{
					\D(\bar x^k,x^k)
				}
			{}\leq{}
				\sum_{k=0}^{K-1}\bigl(
					\FBE(x^k)-\FBE(x^{k+1})
				\bigr)
			{}={}
				\FBE(x^0)-\FBE(x^K)
			{}\leq{}
				\varphi(x^0)-\inf\varphi,
			\end{equation}
			where the last inequality follows from \cref{thm:leq,thm:inf}.
			Since all the iterates up to the \((K-1)\)-th satisfy \(\D(\bar x^k,x^k)>\varepsilon\), if \(\varepsilon>0\) necessarily
			\(
				K
			{}\leq{}
				\frac{\varphi(x^0)-\inf\varphi}{\sigma\varepsilon}
			\)
			as claimed.
		\item\ref{thm:suboptimal}~
			Let
			\(
				\M(z,x)
			{}\coloneqq{}
				f(x)+\innprod{\nabla f(x)}{z-x}
				{}+{}
				g(z)
				{}+{}
				\tfrac1\gamma\D(z,x)
			\)
			be the function in the minimization problem \eqref{eq:T} defining the proximal-gradient operator \(\T\).
			For any \(x\in\mathcal U\) it holds that the difference \(\delta_x(w)\coloneqq\M(w,x)-\varphi(w)\) satisfies
			\[
				\nabla\delta_x(w)
			{}={}
				\nabla(\tfrac1\gamma h-f)(w)
				{}-{}
				\nabla(\tfrac1\gamma h-f)(x).
			\]
			By using convexity of \(L_fh\pm f\) as in the proof of \cref{thm:fC11}, it is easy to verify that the gradient of the (convex) function
			\(
				\tfrac1\gamma h-f
			\)
			is \(\frac{1+\gamma L_f}{\gamma}L_{h,\mathcal U}\)-Lipschitz continuous on \(\mathcal U\), hence so is \(\nabla\delta_x\) independently of \(x\).
			The proof can now trace that of \cite[Lem. 2.15]{themelis2018proximal}.
			Since \(\nabla\delta_x(x)=0\), for any \(\bar x\in\mathcal U\) one has
			\(
				\|\nabla\delta_x(\bar x)\|
			{}\leq{}
				\frac{1+\gamma L_f}{\gamma}L_{h,\mathcal U}\|x-\bar x\|
			\).
			In particular, for \(\bar x\in\T(x)\cap\mathcal U\) one has
			\[
				0
			{}\in{}
				\hat\partial[\M({}\cdot{},x)](\bar x)
			{}={}
				\hat\partial\varphi(\bar x)
				{}+{}
				\nabla\delta_x(\bar x),
			\]
			that is, \(-\nabla\delta_x(\bar x)\in\hat\partial\varphi(\bar x)\).
			Thus, for \(\hat x=\bar x^k\) as in the last iteration of \Cref{alg:Bella} one has
			\[\ifams\mathtight[0.7]\fi
				\dist(0,\hat\partial\varphi(\bar x^k))
			{}\leq{}
				\|\nabla\delta_{x^k}(\bar x^k)\|
			{}\leq{}
				\tfrac{1+\gamma L_f}{\gamma}L_{h,\mathcal U}\|x^k-\bar x^k\|
			{}\leq{}
				\tfrac{1+\gamma L_f}{\gamma}
				\sqrt{
					\tfrac{2L_{h,\mathcal U}^2}{\sigma_{h,\mathcal U}}
					\D(\bar x^k,x^k)
				}
			{}\leq{}
				\tfrac{1+\gamma L_f}{\gamma}
				\sqrt{
					\tfrac{2L_{h,\mathcal U}^2}{\sigma_{h,\mathcal U}}
					\varepsilon
				}
			\]
			as claimed.
			\qedhere
		\end{proofitemize}
	\end{proof}
\end{thm}

We remark that although the proof of \cref{thm:suboptimal} still works even when \(h\) does not have full domain, the result is not very informative in the constrained case \(\overline C\neq\R^n\), owing to the fact that \(0\in\partial\varphi(x_\star)\) is not a necessary condition for optimality when \(x_\star\) is in the boundary of \(\overline C\).
That the result holds regardless is not a contradiction, since the Lipschitz constant \(L_{h,\mathcal U}\) involved in the proof grows bigger and bigger as the boundary is approached.

		\subsection{Subsequential convergence}\label{sec:subseq}
			The rest of the paper is devoted to asymptotic analyses of \cref{alg:Bella}, corresponding to a null tolerance \(\varepsilon=0\).
To rule out trivialities, we will implicitly assume \(\bar x^k\neq x^k\) for every \(k\in\N\), so that the algorithm runs infinitely many iterations.

\begin{thm}[asymptotic analysis]\label{thm:asympt}%
	Suppose that \Cref{ass:basic} holds and consider the iterates generated by \Cref{alg:Bella} with tolerance \(\varepsilon=0\).
	Then,
	\begin{enumerate}
	\item\label{thm:rk}
		\(
			\sum_k\D(\bar x^k,x^k)
		\)
		is finite;
	\item\label{thm:monotone}
		the real-valued sequences \(\seq{\varphi(\bar x^k)}\) and \(\seq{\FBE(x^k)}\) converge to a finite value \(\varphi_\star\), the latter monotonically decreasing;
	\item\label{thm:bounded}
		if \(\varphi+\indicator_C\) is level bounded, then \(\seq{x^k}\) and \(\seq{\bar x^k}\) are bounded;
	\item\label{thm:omega}
		if \(\seq{x^k}\) and \(\seq{\bar x^k}\) are bounded and \(h\) is locally strongly convex, then \(\seq{x^k}\) and \(\seq{\bar x^k}\) have the same set of limit points \(\omega\), which is compact and such that \(\dist(x^k,\omega)\to0\) and \(\dist(\bar x^k,\omega)\to0\) as \(k\to\infty\).
	\end{enumerate}
	\begin{proof}
		\begin{proofitemize}
		\item\ref{thm:rk}~
			Readily follows from the fact that the partial sums in \eqref{eq:telescope} are bounded by the same finite constant for any \(K\in\N\).
		\item\ref{thm:monotone}~
			It follows from \eqref{eq:LS} that \(\seq{\FBE(x^k)}\) is decreasing, hence it admits a limit, be it \(\varphi_\star\), which due to \cref{thm:inf} is lower bounded by \(\inf_C\varphi\) and is thus finite.
			In turn, also \(\varphi(\bar x^k)\to\varphi_\star\) (although not necessarily monotonically), as it follows from \cref{thm:geq} and the fact that \(\D(\bar x^k,x^k)\to0\).
		\item\ref{thm:bounded}~
			Follows from \cref{thm:LB} together with the observation that both sequences are contained both in \(C\) and in the sublevel set \([\FBE\leq\FBE(x^0)]\), as is apparent from \eqref{eq:LS}.
		\item\ref{thm:omega}~
			Let \(\h\coloneqq\tfrac1\gamma h-f\) so that \(\bar x^k\in\prox_\varphi^{\h}\), cf. \cref{thm:Legendre}, and observe that local strong convexity of \(h\) entails that of \(\h\).
			Let \(\sigma_{\h,\mathcal U}\) be a strong convexity modulus for \(\h\) on a (convex) compact set \(\mathcal U\subseteq\dom h\) that contains \(\seq{x^k}\) and \(\seq{\bar x^k}\).
			Then,
			\(
				\tfrac{\sigma_{\h,\mathcal U}}{2}\|x^k-\bar x^k\|^2
			{}\leq{}
				\D*(\bar x^k,x^k)
			\)
			holds for all \(k\), and assertion \ref{thm:rk} then ensures that \(\bar x^k-x^k\to0\).
			We then have
			\[
				x^{k+1}-x^k
			{}={}
				(1-\tau_k)(\bar x^k-x^k)
				{}+{}
				\tau_kd^k
			{}\to{}
				0
			\quad
				\text{as }
				k\to\infty,
			\]
			and since \(x^k-\bar x^k\to0\), similarly
			\(
				\bar x^{k+1}-\bar x^k
			{}\to{}
				0
			\)
			as \(k\to\infty\).
			The claimed properties of \(\omega\) hold generally for any limit set of a bounded sequence.
		\qedhere
		\end{proofitemize}
	\end{proof}
\end{thm}

In the rest of the paper we will need to work under the assumption that \(h\) has full domain.
To the best of our knowledge, no work in the literature dealing with a fully nonconvex setup as that of problem \eqref{eq:P} can circumvent this requirement even for proving stationarity of the limit points.
For the sake of maintaining the full generality of our problem setup we will not include a separate analysis for a (hypo)convex setting and will thus stick to this assumption, which we formulate next for future reference.

\begin{ass}\label{ass:fulldomain}%
	The Legendre function \(h\) has full domain, i.e., \(\dom h=\R^n\).
\end{ass}

\begin{thm}[subsequential convergence]\label{thm:subseq}%
	Suppose that \cref{ass:basic,ass:fulldomain} hold and that \(\varphi\) is level bounded.
	Then, additionally to all the claims of \cref{thm:omega}, any point in the set \(\omega\) satisfies the fixed-point inclusion \(x_\star\in\T(x_\star)\), and in particular is stationary for \(\varphi\).
	Moreover, \(\varphi\) is constant on \(\omega\) and equals \(\varphi_\star\) as in \cref{thm:monotone}.
	\begin{proof}
		It follows from \cref{thm:bounded} that \(\seq{x^k}\) and \(\seq{\bar x^k}\) are bounded, hence that they have cluster points.
		Suppose that \(\seq{\bar x^k}[k\in K]\to x_\star\) for an (infinite) set \(K\subseteq\N\).
		Let \(K'\subseteq K\) be such that \(\bar x^k\to y\) as \(k\in K'\) for some \(y\).
		Then, continuity of \(\D\) and \cref{thm:rk} imply that \(\D(y,x_\star)=0\), hence that \(y=x_\star\).
		From the arbitrarity of the subsequence we conclude that \(\seq{\bar x^k}[k\in K]\to x_\star\).
		By interchanging the roles of \(\bar x^k\) and \(x^k\) it readily follows that \(\seq{x^k}[k\in K]\to x_\star\) iff \(\seq{\bar x^k}[k\in K]\to x_\star\), hence that \(\seq{x^k}\) and \(\seq{\bar x^k}\) have the same set of limit points, be it \(\omega\).
		In particular, any \(x_\star\in\omega\) satisfies \(x_\star\in\T(x_\star)\) owing to \cref{thm:Tosc}.
		Let \(x_\star\in\omega\) be fixed, and let \(K\subseteq\N\) be such that \(x^k,\bar x^k\to x_\star\) as \(K\ni k\to\infty\).
		We have
		\[
			\varphi(x_\star)
		{}\leq{}
			\lim_{K\ni k\to\infty}\varphi(\bar x^k)
		{}\overrel{\ref{thm:monotone}}{}
			\lim_{K\ni k\to\infty}\FBE(x^k)
		{}\overrel{\ref{thm:FBEC0}}{}
			\FBE(x_\star)
		{}\overrel[\leq]{\ref{thm:leq}}{}
			\varphi(x_\star),
		\]
		where the first inequality owes to the fact that \(\varphi\) is lsc.
		The arbitrarity of \(x_\star\in\omega\) and the fact that \(\lim_{k\to\infty}\varphi(\bar x^k)=\varphi_\star\) imply that \(\varphi(x_\star)=\varphi_\star\) for all \(x_\star\in\omega\).
		Let now \(\h\coloneqq\tfrac1\gamma h-f\) so that \(\bar x^k\in\prox_\varphi^{\h}\), cf. \cref{thm:Legendre}.
		The optimality conditions for \(\bar x^k\) read
		\begin{equation}\label{eq:TOC}
			0
		{}\in{}
			\partial\varphi(\bar x^k)
			{}+{}
			\h(\bar x^k)-\nabla\h(x^k).
		\end{equation}
		Continuity of \(\nabla\h\) imply that \(\h(\bar x^k)-\nabla\h(x^k)\to 0\) as \(K\ni k\to\infty\), and since \(\varphi(\bar x^k)\to\varphi(x_\star)\) from \(\varphi\)-attentive osc of \(\partial\varphi\) we conclude that \(0\in\partial\varphi(x_\star)\).
	\end{proof}
\end{thm}

\begin{rem}[adaptive variant of \Cref{alg:Bella} for unknown \(L_f\)]\label{rem:adaptiveLf}%
	If the constant \(L_f\) is not available, then it can be retrieved adaptively by initializing it with an estimate \(L>0\) and by adding the following instruction after \cref{state:BFBS}:
	\begin{algorithmic}[1]
		\algrenewcommand{\alglinenumber}[1]{\footnotesize#1{bis}:}
		\If{
			\(f(\bar x^k) > f(x^k)+\innprod{\nabla f(x^k)}{\bar x^k-x^k}+L\D(\bar x^k,x^k)\)
		}%
			\item[]%
				\(\gamma\gets\nicefrac\gamma2\),~
				\(L\gets 2L\),~
				\(\sigma\gets2\sigma\),~
				and go to \cref{state:BFBS}.%
		\EndIf{}%
	\end{algorithmic}
	Whenever \(L\) exceeds the actual value \(L_f\), this procedure will terminate and \(L\) will be constant starting from that iteration; consequently, \(L\) be increased only a finite number of times.
	Whether or not the final constant \(L\) exceeds the actual value \(L_f\), all the claims of \Cref{thm:subseq} remain valid.
	In order to replicate the proof of \Cref{thm:subseq}, it suffices to show that \(\seq{\FBE(x^k)}\) converges to a finite value \(\varphi_\star\), which here cannot be inferred from the lower boundedness of \(\FBE\) being it ensured only for \(\gamma<\nicefrac{1}{L_f}\) (\Cref{thm:inf}).
	Nevertheless,
	\begin{align*}
		\inf\varphi
	{}\leq{} &
		f(\bar x^k)+g(\bar x^k)
	{}\leq{}
		f(x^k)+\innprod{\nabla f(x^k)}{\bar x^k-x^k}+L\D(\bar x^k,x^k)+g(\bar x^k)
	\\
	{}={} &
		\FBE(x^k)-\tfrac{1-\gamma L}{\gamma}\D(\bar x^k,x^k),
	\end{align*}
	proving that \(\varphi_\star\geq\inf\varphi\).
\end{rem}

		\subsection{Global and linear convergence}\label{sec:global}
			When applying a proximal algorithm to a convex problem, it is typically the case that convergence to a (global) solution is guaranteed whenever one exists.
This appealing property dramatically fails in the nonconvex realm, not only because of the existence of local solutions but because convergence itself (to a single point) may fail even when the iterates are bounded.
A textbook example of this pathological phenomenon is the ``Mexican hat'', a \(\C^\infty\) function defined on the real plane and on which gradient descent may produce sequences with infinitely many limit points \cite{palis1982geometric}.
Luckily, however, under very mild assumptions convergence to a single point can be guaranteed.
The key property is the so-called Kurdyka-Łojasiewicz inequality.

\begin{defin}[KL property]\label{def:KL}%
	A proper lsc function \(\func{F}{\R^n}{\Rinf}\) is said to have the \DEF{Kurdyka-{\L}ojasiewicz} property (KL property) at \(u_\star\in\dom F\) if there exist a concave \DEF{desingularizing function} \(\func{\psi}{[0,\eta]}{[0,\infty)}\) (for some \(\eta>0\)) and an \(\varepsilon>0\) such that
	\begin{enumeratprop}
	\item\label{def:KL1}
		\(\psi(0)=0\);
	\item\label{def:KL2}
		\(\psi\) is of class \(\C^1\) on \((0,\eta)\);
	\item\label{def:KL3}
		for all \(u\in\ball{u_\star}{\varepsilon}\) such that \(F(u_\star)<F(u)<F(u_\star)+\eta\) it holds that
		\begin{equation}\label{eq:KL}
			\psi'(F(u)-F(u_\star))\dist(0,\partial F(u))\geq 1.
		\end{equation}
	\end{enumeratprop}
\end{defin}

The first inequality of this type is given in the seminal work of {\L}ojasiewicz  \cite{lojasiewicz1963propriete,lojasiewicz1993geometrie} for analytic functions, which we nowadays call {\L}ojasiewicz's gradient inequality.
Kurdyka \cite{kurdyka1998gradients} showed that this inequality is valid for \(\C^1\) functions whose graph belongs to an \DEF{\(o\)-minimal structure} \cite{van1996geometric,van1998tame}, a result which has later been extended in \cite{bolte2007lojasiewicz,bolte2007clarke,bolte2010characterizations} for lsc \DEF{tame} functions, a wide category including semialgebraic functions as special case.
In the following result we adapt the convergence analysis of \cite{attouch2010proximal,attouch2013convergence,noll2014convergence} to our framework, showing that \cref{alg:Bella} converges to a (unique) limit point under a KL property assumption.

\begin{thm}[global convergence]\label{thm:global}%
	Let \Cref{ass:basic,ass:fulldomain} hold and consider the iterates generated by \Cref{alg:Bella} with tolerance \(\varepsilon=0\).
	Suppose further that
	\begin{enumeratass}
	\item
		\(\varphi\) is level bounded;
	\item
		\(f,h\in\C^2\) with \(\nabla^2h\succ0\);
	\item\label{ass:dk}
		either \(\seq{\|d^k\|}\) has finite sum, or there exists \(D\geq 0\) such that \(\|d^k\|\leq D\|x^k-\bar x^k\|\) for all \(k\);
	\item
		\(f,g,h\) are tame functions \cite{van1996geometric,van1998tame} (e.g., semialgebraic).
	\end{enumeratass}
	Then, \(\seq{\|x^k-\bar x^k\|}\) has finite sum (in fact, regardless of whether or not \cref{ass:dk} holds), and \(\seq{x^k}\) and \(\seq{\bar x^k}\) converge to a stationary point \(x_\star\).
	\begin{proof}
		Let \(\h\coloneqq\frac1\gamma h-f\) so that
		\(
			\bar x^k
		{}={}
			\prox_\varphi^{\h}(x^k)
		{}={}
			\argmin_w\M(w,x)
		\),
		where \(\M(w,x)\coloneqq\varphi(w)+\D*(w,x)\) (cf. \cref{thm:Legendre}).
		It follows from \cref{thm:monotone} that \(\FBE(x^k)=\M(\bar x^k,x^k)\) converges strictly decreasing to \(\varphi_\star\).
		Since tame functions are closed under derivation and basic algebraic operations \cite[\S2.1]{van1996geometric}, \(\M\) is
		\ifams
			tame
		\else
			a tame function
		\fi
		and consequently satisfies the KL property \cite[Thm. 14]{bolte2007clarke}.
		Moreover, \(\M\) is constant (and equals \(\varphi_\star\)) on \(\Omega\coloneqq\set{(x_\star,x_\star)}[x_\star\in\omega]\).
		The properties of the set of accumulation points \(\omega\) of the sequences \(\seq{x^k}\) and \(\seq{\bar x^k}\) asserted in \cref{thm:omega} ensure through \cite[Lem. 6]{bolte2014proximal} the existence of a \emph{uniformized KL function} on \(\Omega\), namely a function \(\psi\) satisfying \cref{{def:KL1},,{def:KL2},,{def:KL3}} with \(F=\M\) for all \(u_\star\in\Omega\) and \(u\in\R^n\times\R^n\) such that
		\(
			\dist(u,\Omega)<\varepsilon
		\)
		and
		\(
			\varphi_\star < \M(u) < \varphi_\star + \eta
		\).
		Up to possibly discarding the first iterates, we may assume that \(u=(\bar x^k,x^k)\) satisfies these conditions for all \(k\in\N\).
		Let \(\Delta_k\coloneqq\psi(\M(\bar x^k,x^k)-\varphi_\star)=\psi(\FBE(x^k)-\varphi_\star)\), and observe that
		\begin{equation}\label{eq:partialM}
			\partial\M(w,x)
		{}={}
			\binom{
				\partial\varphi(w)+\nabla\h(w)-\nabla\h(x)
			}{
				\nabla^2\h(x)(x-w)
			},
		\quad\text{hence}\quad
			\binom{0}{\nabla^2\h(x^k)(x^k-\bar x^k)}
		{}\in{}
			\partial\M(\bar x^k,x^k)
		\end{equation}
		as it follows from \eqref{eq:TOC}.
		We have
		\begin{align*}
			1
		{}\leq{} &
			\psi'(\M(\bar x^k,x^k)-\varphi_\star)
			\dist(0,\partial\M(\bar x^k,x^k))
		&&
			\dueto{KL inequality \eqref{eq:KL}}
		\\
		\numberthis\label{eq:MKL}
		{}\leq{} &
			\psi'(\M(\bar x^k,x^k)-\varphi_\star)
			\|\nabla^2\h(x^k)\|\|x^k-\bar x^k\|
		&&
			\dueto{inclusion \eqref{eq:partialM}}
		\\
		{}\leq{} &
			\frac{\Delta_k-\Delta_{k+1}}{\M(\bar x^k,x^k)-\M(\bar x^{k+1},x^{k+1})}
			\|\nabla^2\h(x^k)\|\|x^k-\bar x^k\|
		&&
			\dueto{concavity of \(\psi\)}
		\\
		{}\leq{} &
			\frac{\Delta_k-\Delta_{k+1}}{\sigma\D(\bar x^k,x^k)}
			\|\nabla^2\h(x^k)\|\|x^k-\bar x^k\|
		&&
			\dueto{\(\M(\bar x^j,x^j)=\FBE(x^j)\) and \eqref{eq:LS}}
		\\
		\numberthis\label{eq:Deltadiff}
		{}\leq{} &
			\frac{2L_{\h,\mathcal U}}{\sigma\sigma_{\h,\mathcal U}}
			\frac{\Delta_k-\Delta_{k+1}}{\|x^k-\bar x^k\|},
		\end{align*}
		where \(L_{\h,\mathcal U}\) and \(\sigma_{\h,\mathcal U}\) are, respectively, smoothness and strong convexity moduli of \(\h\) on a (convex) compact set \(\mathcal U\) that contains all the iterates \(x^k\) and \(\bar x^k\).
		Therefore,
		\begin{equation}\label{eq:rksummable}
			\sum_{k\in\N}\|x^k-\bar x^k\|
		{}\leq{}
			\frac{2L_{\h,\mathcal U}}{\sigma\sigma_{\h,\mathcal U}}
			\sum_{k\in\N}(\Delta_k-\Delta_{k+1})
		{}\leq{}
			\frac{2L_{\h,\mathcal U}}{\sigma\sigma_{\h,\mathcal U}}
			\Delta_0
		{}<{}
			\infty,
		\end{equation}
		where the last inequality follows from the fact that \(\psi\geq0\).
		Therefore,
		\[
			\sum_{k\in\N}\|x^{k+1}-x^k\|
		{}={}
			\sum_{k\in\N}\|(1-\tau_k)(\bar x^k-x^k)+\tau_kd^k\|
		{}\leq{}
			\sum_{k\in\N}\|\bar x^k-x^k\|
			{}+{}
			\sum_{k\in\N}\|d^k\|
		{}<{}
			\infty
		\]
		as ensured by either conditions in \cref{ass:dk} in light of \cref{eq:rksummable}.
		In particular, \(\seq{x^k}\) is a Cauchy sequence and thus has a limit \(x_\star\), which is also the limit of \(\seq{\bar x^k}\) and is stationary for \(\varphi\) as it follows from \cref{thm:subseq}.
	\end{proof}
\end{thm}

\Cref{ass:dk} imposes a mild and reasonable consistency criterion on the directions \(d^k\) being used (which however is not needed for subsequential convergence, cf. \cref{thm:subseq}).
It simply reflects the idea that shorter steps should be taken when close to solutions, and uses the fixed-point residual \(\|x^k-\bar x^k\|\) to quantify the proximity.
Though apparently quite abstract a requirement, in the next subsection a more practical understanding of this condition will be given when showing its connection to standard theory of quasi-Newton schemes.

When \(f,g,h\) are semialgebraic (and thus so is the model \(\M\) in the proof of \cref{thm:global} as a byproduct), then the desingularizing function for \(\M\) can be taken of the form \(\psi(s)=\varrho s^\vartheta\) for some \(\varrho>0\) and \(\vartheta\in(0,1)\) \cite{attouch2010proximal}, in which case we say that it satisfies the KL property \emph{with exponent \(1-\vartheta\)}.
Such exponent has a cardinal role in determining asymptotic rates of convergence; in particular, when \(\vartheta\geq\nicefrac12\) linear convergence rates can easily be shown.
We remark that our approach based on the model \(\M\) (both in \cref{thm:global} and in \cref{thm:linear} that follows) is largely inspired by the similar approach in \cite{li2016douglas} for the global and linear convergence of the Douglas-Rachford algorithm in the nonconvex case.

\begin{thm}[{linear convergence}]\label{thm:linear}%
	Additionally to the assumptions of \Cref{thm:global}, suppose that \(f,g,h\) are semialgebraic, and that the KL exponent of \(\M(w,x)\coloneqq\varphi(w)+\D*(w,x)\) with \(\h\coloneqq\tfrac1\gamma h-f\) is \(1-\vartheta\leq\nicefrac12\).
	Then, \(\seq{x^k}\) and \(\seq{\bar x^k}\) generated by \cref{alg:Bella} converge at \(R\)-linear rate to a stationary point.
	\begin{proof}
		As shown in \Cref{thm:global}, the sequences converge to stationary point \(x_\star\).
		Since \(x^{k+1}-x^k=(1-\tau_k)(\bar x^k-x^k)+\tau_kd^k\), defining \(B_k\coloneqq\sum_{i\geq k}\|\bar x^i-x^i\|\) one has
		\(
			\|x^k-x_\star\|
		{}\leq{}
			(1+c)B_k
		\),
		and similarly
		\(
			\|\bar x^k-x_\star\|
		{}\leq{}
			(3+c)B_k
		\)
		owing to the inequality
		\[
			\|\bar x^{k+1}-\bar x^k\|
		{}\leq{}
			\|\bar x^{k+1}-x^{k+1}\|
			{}+{}
			\|\bar x^k-x^k\|
			{}+{}
			\|x^{k+1}-x^k\|
		{}\leq{}
			\|\bar x^{k+1}-x^{k+1}\|
			{}+{}
			(2+c)\|\bar x^k-x^k\|.
		\]
		As such, it suffices to show that the sequence \(\seq{B_k}\) converges with asymptotic \(R\)-linear rate.
		Let \(\Delta_k\), \(L_{\h,\mathcal U}\) and \(\sigma_{\h,\mathcal U}\) be as in the proof of \cref{thm:global}.
		The KL inequality \eqref{eq:MKL} with \(\psi(s)=\varrho s^\vartheta\) reads
		\[
			1
		{}\leq{}
			\varrho\vartheta L_{\h,\mathcal U}
			(\M(\bar x^k,x^k)-\varphi_\star)^{\vartheta-1}
			\|x^k-\bar x^k\|
		{}={}
			\varrho\vartheta L_{\h,\mathcal U}
			(\FBE(x^k)-\varphi_\star)^{\vartheta-1}
			\|x^k-\bar x^k\|.
		\]
		Since \(\|x^k-\bar x^k\|\to0\), up to discarding the first iterates we may assume that this quantity is smaller than 1.
		Therefore, \(\Delta_k=\psi(\FBE(x^k)-\varphi_\star)\) satisfies
		\begin{equation}\label{eq:Bkrk}
			\Delta_k
		{}={}
			\varrho(\FBE(x^k)-\varphi_\star)^\vartheta
		{}\leq{}
			\varrho(\varrho\vartheta L_{\h,\mathcal U})^{\frac{\vartheta}{1-\vartheta}}
			\|x^k-\bar x^k\|^{\frac{\vartheta}{1-\vartheta}}
		{}\leq{}
			\varrho^{\frac{1}{1-\vartheta}}
			(\vartheta L_{\h,\mathcal U})^{\frac{\vartheta}{1-\vartheta}}
			\|x^k-\bar x^k\|,
		\end{equation}
		where the last inequality uses the fact \(\frac{\vartheta}{1-\vartheta}\geq1\).
		Hence,
		\[
			B_k
		{}={}
			\sum_{i\geq k}{
				\|x^i-\bar x^i\|
			}
		{}\overrel*[\leq]{\eqref{eq:Deltadiff}}{}
			\frac{2L_{\h,\mathcal U}}{\sigma\sigma_{\h,\mathcal U}}
			\sum_{i\geq k}(
				\Delta_i-\Delta_{i+1}
			)
		{}\overrel*[\leq]{\(\mathtight\Delta_i\geq0\)}{}
			\frac{2L_{\h,\mathcal U}}{\sigma\sigma_{\h,\mathcal U}}
			\Delta_k
		{}\overrel*[\leq]{\eqref{eq:Bkrk}}{}
			c\|x^k-\bar x^k\|
		\]
		for some constant \(c>0\).
		Therefore,
		\(
			B_k
		{}\leq{}
			c\|x^k-\bar x^k\|
		{}={}
			c(B_k-B_{k+1})
		\),
		leading to the sought linear rate
		\(
			B_{k+1}
		{}\leq{}
			(1-\nicefrac1c)
			B_k
		\).
	\end{proof}
\end{thm}

		\subsection{Superlinear convergence}\label{sec:superlinear}
			Although \cref{alg:Bella} is ``robust'' to any choice of directions, a suitable selection stemming for instance from Newton-type method can cause a remarkable speed-up.
As already discussed in \cref{thm:2nd}, the forward-backward mapping \(\T\) behaves nicely around nondegenerate fixed points under some regularity assumptions.
This motivates the quest to derive the direction \(d^k\) in \Cref{alg:Bella} by a Newton-type scheme for solving the inclusion \(x\in\T(x)\) or, equivalently, the nonlinear (generalized) equation \(0\in \Res(x)\) where \(\Res\coloneqq\id-\T\) is the fixed-point residual mapping.
Newton-type methods for such a nonlinear equation prescribe updates of the form
\begin{equation}\label{eq:Newton}
    x^+=x-H(x)\Res(x),
\end{equation}
where ideally \(H(x)\) should well approximate the Jacobian of \(\Res(x)\).
In particular, starting with an invertible matrix \(H_0\), quasi-Newton schemes emulate higher-order information by performing low-rank updates satisfying a so-called secant equation
\begin{equation}\label{eq:Secant}
    H^+y=s,
\quad
	\text{where }
	s=x^+-x,~
	y\in\Res(x^+)-\Res(x).
\end{equation}
In these scenarios, the condition \(\|d^k\|\leq D\|x^k-\bar x^k\|\) appearing in \cref{ass:dk} is verified when the operators \(H\) are bounded, a frequent assumption in the convergence analysis of quasi-Newton methods.
A well-known result characterizing the superlinear convergence of this type of schemes is based on the Dennis-Mor\'e condition \cite{dennis1974characterization,dennis1977quasi}, which amounts to differentiability of \(\Res\) at the limit point together with the limit \(\|\Res(x^k)+J\Res(x_\star)d^k\|/\|d^k\|\to0\); see also \cite{dontchev2012generalizations} for the extension to generalized equations.
In \Cref{thm:DM}, we will see that directions satisfying this condition do trigger asymptotic superlinear rates in \Cref{alg:Bella}.
To this end, we first characterize the quality of the update directions with the next definition, and prove an intermediate result showing how they fit into \Cref{alg:Bella}.
In the sequel, we will make use of the notion of \DEF{nonisolated superlinear directions} that we define next.

\begin{defin}[nonisolated superlinear directions]\label{def:supDir}%
	Relative to the iterates generated by \Cref{alg:Bella}, we say that \(\seq{d^k}\) is a sequence of \DEF{superlinear directions with order \(q\geq1\)} if
	\[
		\lim_{k\to\infty}{
			\frac{\dist(x^k+d^k,\mathcal X_\star)}{\dist(x^k,\mathcal X_\star)^q}
		}
	{}={}
		0,
	\]
	with
	\(\mathtight[0.6]
		\mathcal X_\star
	{}\coloneqq{}
		\set{x\in\R^n}[x\in\T(x)]
	\)
	the set of fixed points of the forward-backward operator \(\T\).%
\end{defin}

The set \(\mathcal X_\star\) in the definition above corresponds to the possible limit points of the sequences \(\seq{x^k}\) and \(\seq{\bar x^k}\) generated by \cref{alg:Bella} when \(h\) has full domain, as shown in \cref{thm:subseq}.
In fact, our notion of superlinear directions extends the one given in \cite[\S7.5]{facchinei2003finite} to cases in which \(\mathcal X_\star\) is not a singleton.

\subsubsection{A nonlinear error bound}
	Our result hinges on three key ingredients: (1) the differentiability of the BFBE around a limit point of \cref{alg:Bella}, (2) the KL property (with exponent) on the BFBE, and (3) a nonlinear error bound relating the KL function to the distance appearing in \cref{def:supDir}.
	Sufficient conditions ensuring the first ingredient have already been discussed in \cref{thm:FBEC1}.
	As to the second one, it has been shown in \cite{yu2019deducing} that whenever the kernel function \(h\) is twice continuously differentiable and (locally) strongly convex, if the function \(\varphi\) satisfies the KL property with exponent \(\vartheta\geq\nicefrac12\) then so does the Bregman envelope \(\varphi^h\).
	Clearly, the lower the \(\vartheta\) the stronger the property, in the sense that whenever \(\varphi\) admits a desingularizing function with exponent \(\vartheta\in(0,1)\), then it also admits a desingularizing function with exponent \(\vartheta'\in[\vartheta,1)\).
	Combined with the relation existing among the BFBE and the Bregman-Moreau envelope shown in \cref{thm:Legendre}, the following is obtained.

	\begin{lem}[equivalence of \L ojasiewicz property {\cite[Thm. 5.2]{yu2019deducing}}]\label{thm:phiKL}%
		Additionally to \Cref{ass:basic}, suppose that \(f,h\in\C^2(\interior\dom h)\) with \(\nabla^2h\succ 0\).
		If \(\varphi\) has the KL property with exponent \(\vartheta\in(0,1)\), then so does \(\FBE\) with exponent \(\vartheta'=\max\set{\nicefrac12,\vartheta}\).
	\end{lem}

	The last ingredient involves a desingularizing property stronger than the KL inequality which is investigated in \cite{aze2017nonlinear}, namely with \(\dist(0,\partial F(u))\) being replaced by the \DEF{strong slope}
	\(
		|\nabla F|(i)
	{}\coloneqq{}
		\limsup_{u\neq z\to u}\frac{F(u)-F(z)}{\|u-z\|}
	\)
	in \eqref{eq:KL}.
	Of particular interest to our scope is \cite[Thm. 4.1]{aze2017nonlinear}, which connects this stronger property to a nonlinear growth condition of the form \(\psi(F(u)-F(u_\star))\geq\dist(u,[F\leq F(u_\star)])\).
	The key observation is that whenever \(F\) is continuously differentiable around \(u_\star\), both the strong slope \(|\nabla F|(u)\) and the minimum norm subgradient \(\dist(0,\partial F(u))\) coincide and are equal to \(\|\nabla F(u)\|\).
	The combination of \cite[Thm. 4.1]{aze2017nonlinear} with \cref{thm:FBEC1} thus leads to the following milestone of our analysis.

	\begin{lem}[nonlinear error bound {\cite[Thm. 4.1]{aze2017nonlinear}}]\label{thm:EB}%
		Suppose that \Cref{ass:basic} holds and let \(x_\star\) be a nondegenerate fixed point of \(\T\) at which \(g\) is prox-regular for \(-\nabla f(x_\star)\).
		Suppose further that \(f,h\in\C^2(\interior\dom h)\) with \(\nabla^2h\succ0\), and that \(\FBE\) has the KL property at \(x_\star\) with desingularizing function \(\psi\).
		Then, denoting \(\varphi_\star\coloneqq\FBE(x_\star)\) there exist \(\varepsilon,\eta>0\) such that
		\[
			\psi\bigl(\FBE(x)-\varphi_\star\bigr)
		{}\geq{}
			\dist\bigl(x,[\FBE\leq\varphi_\star]\bigr)
		\quad
			\forall x\in\ball{x_\star}{\varepsilon}
			\text{ such that }
			\varphi_\star<\FBE(x)<\varphi_\star+\eta.
		\]
	\end{lem}

\subsubsection{Superlinear convergence to nonisolated local minima}
	We now have all the ingredients to address the superlinear convergence analysis of \cref{alg:Bella} to nonisolated local minima.
	The next result constitutes a key component of the proposed methodology, as it shows that \cref{alg:Bella} does not suffer from the \emph{Maratos' effect} \cite{maratos1978exact}, a well-known obstacle for fast local methods that inhibits the acceptance of the unit stepsize.
	On the contrary, we will show that under mild assumptions whenever the directions \(\seq{d^k}\) in \Cref{alg:Bella} are superlinear, then indeed unit stepsize is eventually always accepted and the algorithm converges superlinearly even if the limit point belongs to a flat region of local minima.
	As detailed in the proof, local minimality (as opposed to the more general property of being a fixed point of the forward-backward mapping \(\T\)) enables the identity \(\dist\bigl(x,[\FBE\leq\varphi_\star]\bigr)=\dist(x,\mathcal X_\star)\) in a neighborhood of the limit point, which in turns allows us to connect the notion of superlinear directions of \cref{def:supDir} to the nonlinear error bound of \cref{thm:EB}.
	
	\begin{thm}[acceptance of the unit stepsize and superlinear convergence]\label{thm:tau1}%
		Consider the iterates generated by \Cref{alg:Bella}, and additionally to \Cref{ass:basic,ass:fulldomain} suppose that the following requirements hold:
		\begin{enumeratass}
		\item
			\(\varphi\) is level bounded;
		\item
			\(f,h\in\C^2(\R^n)\) with \(\nabla^2h\succ0\);
		\item
			\(\varphi\) has the KL property with exponent \(1-\vartheta\in(0,1)\) (as is the case when \(f,g,h\) are semialgebraic);
		\item
			\(d^k\) are superlinear directions with order \(q\geq\max\set{1,\nicefrac{1}{2\vartheta}}\) (cf. \cref{def:supDir});
		\item
			the sequence \(\seq{x^k}\) converges to a nondegenerate fixed point \(x_\star\) of \(\T\), which is a (not necessarily isolated) local minimum for \(\varphi\), and at which \(g\) is prox-regular for \(-\nabla f(x_\star)\).
		\end{enumeratass}
		Then, there exists \(k_0\in\N\) such that
		\[
			\FBE(x^k+d^k)
		{}\leq{}
			\FBE(x^k)-\sigma\D(\bar x^k,x^k)
		\quad
			\forall 
			k\geq k_0.
		\]
		In particular,
		\begin{enumerate}
		\item
			eventually stepsize \(\tau=1\) is always accepted at \cref{state:LS} (that is, no backtrackings eventually occur) and the iterates reduce to \(x^{k+1}=x^k+d^k\);
		\item
			\(\dist(x^k,\mathcal X_\star)\to0\) at superlinear rate, where \(\mathcal X_\star\) is as in \Cref{def:supDir}.
		\end{enumerate}
		\begin{proof}
			Firstly, \cref{thm:phiKL} ensures that \(\psi(s)\coloneqq\varrho s^{\min\set{\vartheta,\nicefrac12}}\) for some \(\varrho>0\) is a desingularizing function for \(\FBE\) at \(x_\star\).
			Denoting \(\varphi_\star\coloneqq\varphi(x_\star)=\FBE(x_\star)\), the equivalence of local minimality asserted in \cref{thm:locmin} ensures that for small enough \(\varepsilon>0\) it holds that
			\begin{equation}\label{eq:lev=1}
				\FBE(x)
			{}\geq{}
				\varphi_\star
			~~\text{and}~~
				\dist\bigl(x,[\FBE\leq\varphi_\star]\bigr)
			{}={}
				\dist\bigl(x,[\FBE=\varphi_\star]\bigr)
			~~\forall
				x\in\ball{x_\star}{\varepsilon}.
			\end{equation}
			Moreover, since \(\FBE(x^k)\) converges strictly decreasing to \(\varphi_\star\) (cf. \cref{thm:monotone}), up to possibly discarding the first iterates and restricting \(\varepsilon\) we may assume that \(\FBE(x^k)\leq\varphi_\star+\eta\), with \(\eta,\varepsilon>0\) as in \cref{def:KL} of the KL function.
			Notice further that any point \(x\in\ball{x_\star}{\varepsilon}\) and such that \(\FBE(x)=\varphi_\star\) necessarily satisfies \(x\in\mathcal X_\star\) owing to \cref{thm:leq} and local minimality of \(x_\star\).
			Conversely, it follows from \cref{thm:FBEC1} that \(\FBE\) is \(\C^1\) around \(x_\star\) with \(\nabla\FBE(x)=\Q(x)(x-\T(x))\), where \(\Q\succ0\) is as in \eqref{eq:Q}, and in particular \(\nabla\FBE(x)=0\) for any \(x\in\mathcal X_\star\) close to \(x_\star\).
			Combined with the KL inequality \eqref{eq:KL}, we conclude that close to \(x_\star\) a point \(x\) belongs to \(\mathcal X_\star\) iff \(\varphi(x)=\FBE(x)=\varphi_\star\).
			Up to possibly further restricting \(\varepsilon\), we may thus modify \eqref{eq:lev=1} to
			\begin{equation}\label{eq:lev=}
				\dist\bigl(x,[\FBE=\varphi_\star]\bigr)
			{}={}
				\dist\bigl(x,\mathcal X_\star\bigr)
			\quad
				\forall
				x\in\cball{x_\star}{\varepsilon}.
			\end{equation}
			Combined with the error bound in \cref{thm:EB} with \(\psi(s)=\varrho s^{\min\set{\vartheta,\nicefrac12}}\), we obtain
			\begin{equation}\label{eq:phidist}
				\FBE(x^k)-\varphi_\star
			{}\geq{}
				\left(\varrho^{-1}\dist(x^k,\mathcal X_\star)\right)^{\max\set{2,\nicefrac1\vartheta}}
			\quad
				\forall k\in\N.
			\end{equation}
			Since, as discussed above, \(\mathcal X_\star\) coincides with a (closed) sublevel set of \(\FBE\) close to \(x_\star\), for every \(k\) there exists a projection point \(x_\star^k\in\proj_{\mathcal X_\star}(x^k+d^k)\), hence such that \(\varphi(x_\star^k)=\varphi_\star\) as motivated before.
			In particular,
			\begin{align*}
				\FBE(x^k+d^k)
			{}\leq{} &
				\varphi(x_k^\star)+\D*(x_k^\star,x^k+d^k)
			{}\leq{}
				\varphi_\star
				{}+{}
				\tfrac{L_{\h,\mathcal U}}{2}
				\dist(x^k+d^k,\mathcal X_\star)^2,
			\end{align*}
			where \(\h\coloneqq\tfrac1\gamma h-f\in\C^2(\R^n)\), \(L_{\h,\mathcal U}\) is a Lipschitz modulus of \(\nabla\h\) on \(\mathcal U=\cball{x_\star}{\varepsilon}\), and the first inequality follows from \cref{thm:Legendre}.
			Together with \eqref{eq:phidist}, this implies
			\begin{equation}\label{eq:rkLim}
				\varepsilon_k
			{}\coloneqq{}
				\frac{
					\FBE(x^k+d^k)-\varphi_\star
				}{
					\FBE(x^k)-\varphi_\star
				}
			{}\leq{}
				\frac{L_{\h,\mathcal U}}{2}
				\left(
					\frac{
						\dist(x^k+d^k,\mathcal X_\star)
					}{
						\left(\varrho^{-\nicefrac12}\dist(x^k,\mathcal X_\star)\right)^{\max\set{1,\nicefrac{1}{2\vartheta}}}
					}
				\right)^2
			{}\to{}
				0
			\quad\text{as }
				k\to\infty.
			\end{equation}
			Thus, for large enough \(k\) so that \(\varepsilon_k\leq1\), we have
			\begin{align*}
				\FBE(x^k+d^k)
				{}-{}
				\FBE(x^k)
			{}={} &
				\bigl(
				\FBE(x^k+d^k)
					{}-{}
					\varphi_\star
				\bigr)
				{}-{}
				\bigl(
					\FBE(x^k)
					{}-{}
					\varphi_\star
				\bigr)
			\\
			{}={} &
				(\varepsilon_k-1)
				\bigl(
					\FBE(x^k)
					{}-{}
					\varphi_\star
				\bigr)
			\\
				\dueto{(since
					\(\varepsilon_k\leq1\),~
					\(\varphi(\bar x^k)\geq\varphi_\star\)%
				)}~
			{}\leq{} &
				(\varepsilon_k-1)
				\bigl(
					\FBE(x^k)
					{}-{}
					\varphi(\bar x^k)
				\bigr)
			\\
				\dueto{(use \cref{thm:geq})}~
			{}\leq{} &
				(\varepsilon_k-1)
				\tfrac{1-\gamma L_f}{\gamma}\D(\bar x^k,x^k)
			{}\leq{}
				\sigma\D(\bar x^k,x^k)
			\end{align*}
			for large enough \(k\), where the last inequality holds since \(\sigma<\tfrac{1-\gamma L_f}{\gamma}\) and \(\varepsilon_k\to 0\).
		\end{proof}
	\end{thm}

\subsubsection{Superlinear convergence to strong minima under the Dennis-Moré condition}
	Despite the importance of nonisolated critical points in nonsmooth nonconvex optimization, there has been little about superlinear directions for such problems. In the convex setting, some studies have shown the potential of variants of regularized Newton \cite{li2004regularized,themelis2019acceleration} and semismooth Newton methods \cite{li2018highly,themelis2019supermann} under a local error bound.
	In the smooth nonconvex setting, there are many works relying on Levenberg-Marquardt \cite{ahookhosh2019local,ahookhosh2020finding,fischer2002local,yamashita2001rate}, cubic regularization \cite{yue2019quadratic}, and regularized Newton \cite{ueda2014regularized} methods under variants of local error bounds and Hölder metric subregularity.

	Quasi-Newton methods constitute an important class of directions widely used in optimization.
	Superlinear convergence of these type of directions is typically assessed by means of the Dennis-Moré condition.
	We next show that under regularity assumptions at the limit point the same condition ensures acceptance of unit stepsize in our framework, albeit provided the algorithm converges to an (isolated) strong local minimum.
	%

	\begin{thm}[superlinear convergence under Dennis-Moré condition]\label{thm:DM}%
		Consider the iterates generated by \Cref{alg:Bella}.
		Additionally to \Cref{ass:basic,ass:fulldomain}, suppose that the following requirements are satisfied:
		\begin{enumeratass}
		\item
			\(\seq{x^k}\) converges to a strong local minimum \(x_\star\) of \(\varphi\);
		\item
			\(f,h\in\C^2(\R^n)\) with \(\nabla^2h\succ0\);
		\item
			\(\Res(x)\coloneqq x-\T(x)\) is strictly differentiable at \(x_\star\) (see \cref{thm:2nd} for sufficient conditions) with \(J\Res(x_\star)\) nonsingular;
		\item
			\(\seq{d^k}\) satisfy the Dennis-Moré condition
			\begin{equation}\label{eq:DM}
				\lim_{k\to\infty}{
					\frac{
						\Res(x^k)+J\Res(x_\star)d^k
					}{
						\|d^k\|
					}
				}
			{}={}
				0.
			\end{equation}
		\end{enumeratass}
		Then, \(\seq{d^k}\) are superlinear directions with order \(q=1\), and in particular all the claims of \Cref{thm:tau1} hold.
		\begin{proof}
			Since \(\Res(x_\star)=0\) as it follows from \cref{thm:subseq}, the Dennis-Moré condition \eqref{eq:DM} and strict differentiability at \(x_\star\) imply that
			\[
				\lim_{k\to\infty}{
					\frac{
						\Res(x^k+d^k)
					}{
						\|d^k\|
					}
				}
			{}={}
				\lim_{k\to\infty}\left[
					\frac{
						\Res(x^k)+J\Res(x_\star)d^k-\Res(x^k+d^k)
					}{
						\|d^k\|
					}
					{}+{}
					\frac{
						\Res(x^k+d^k)
					}{
						\|d^k\|
					}
				\right]
			{}={}
				0.
			\]
			Further, nonsingularity of \(J\Res(x_\star)\) implies that there exists \(\alpha>0\) such that
			\[
				\|\Res(x)\|
			{}={}
				\|\Res(x)-\Res(x_\star)\|
			{}\geq{}
				\alpha\|x-x_\star\|
			\]
			holds for all \(x\) close enough to \(x_\star\).
			We thus have
			\[
				\ifams
					\mathtight[0.8]
				\fi
				0
			{}\leftarrow{}
				\frac{
					\|
						\Res(x^k+d^k)
					\|
				}{
					\|d^k\|
				}
			{}\geq{}
				\alpha\frac{
					\|
						x^k+d^k-x_\star
					\|
				}{
					\|d^k\|
				}
			{}\geq{}
				\alpha\frac{
					\|
						x^k+d^k-x_\star
					\|
				}{
					\|x^k+d^k-x_\star\|
					{}+{}
					\|x^k-x_\star\|
				}
			{}={}
				\alpha\frac{
					\frac{
						\|
							x^k+d^k-x_\star
						\|
					}{
						\|x^k-x_\star\|
					}
				}{
					1
					{}+{}
					\frac{
						\|
							x^k+d^k-x_\star
						\|
					}{
						\|x^k-x_\star\|
					}
				},
			\]
			as \(k\to\infty\), and in particular
			\(
				\frac{
					\|
						x^k+d^k-x_\star
					\|
				}{
					\|x^k-x_\star\|
				}
			{}\to{}
				0
			\),
			as claimed.
		\end{proof}
	\end{thm}

 	\section{Final remarks}\label{sec:Conclusion}

We proposed \Cref{alg:Bella}, a Bregman-forward-backward-splitting-based algorithm for minimizing the sum of two nonconvex functions, where the first one is relatively smooth and the second one is possibly nonsmooth.
\Cref{alg:Bella} is a linesearch algorithm on the Bregman forward-backward envelope (BFBE), a Bregman extension of the forward-backward envelope, and globalizes convergence of fast local methods for finding zeros of the forward-backward residual.
Furthermore, thanks to a nonlinear local error bound holding for the BFBE under prox-regularity and the KL property, the algorithm enables acceptance of unit stepsize when the fast local method yields directions that are superlinear with respect to the set of solutions, thus triggering superlinear convergence even when the limit point is not isolated.

In future work we plan to address the following issues:
(1) extending existing superlinear direction schemes such as those proposed in \cite{li2004regularized,ueda2014regularized,ahookhosh2019local,themelis2019acceleration} for either convex or smooth problems to the more general setting of this paper;
(2) assessing the performance of such schemes in the \Cref{alg:Bella} framework with numerical simulations on nonconvex nonsmooth problems such as low-rank matrix completion, sparse nonnegative matrix factorization, phase retrieval, and deep learning; and
(3) guaranteeing saddle point avoidance, in the spirit of \cite{oneill2018behavior,lee2019first,liu2019envelope}.


\ifsiam
	\bibliographystyle{siamplain}
\else
	\bibliographystyle{plain}
\fi
	\bibliography{TeX/Bibliography.bib}

\begin{thebibliography}{10}

\bibitem{ahookhosh2018optimal}
Masoud Ahookhosh.
\newblock Optimal subgradient methods: computational properties for large-scale
  linear inverse problems.
\newblock {\em Optimization and Engineering}, 19(4):815--844, 2018.

\bibitem{ahookhosh2019accelerated}
Masoud Ahookhosh.
\newblock Accelerated first-order methods for large-scale convex optimization:
  nearly optimal complexity under strong convexity.
\newblock {\em Mathematical Methods of Operations Research}, 89(3):319--353,
  2019.

\bibitem{ahookhosh2019local}
Masoud Ahookhosh, Francisco J~Arag\'on Artacho, Ronan~MT Fleming, and Phan~T
  Vuong.
\newblock Local convergence of the {L}evenberg--{M}arquardt method under
  {H}{\"o}lder metric subregularity.
\newblock {\em Advances in Computational Mathematics}, pages 1--36, 2019.

\bibitem{ahookhosh2020finding}
Masoud Ahookhosh, Ronan~MT Fleming, and Phan~T Vuong.
\newblock Finding zeros of {H}{\"o}lder metrically subregular mappings via
  globally convergent {L}evenberg--{M}arquardt methods.
\newblock {\em Optimization Methods and Software}, pages 1--37, 2020.

\bibitem{ahookhosh2019multi}
Masoud Ahookhosh, Le~Thi~Khanh Hien, Nicolas Gillis, and Panagiotis Patrinos.
\newblock Multi-block {B}regman proximal alternating linearized minimization
  and its application to orthogonal nonnegative matrix factorization.
\newblock {\em arXiv:1908.01402}, 2019.

\bibitem{ahookhosh2020block}
Masoud Ahookhosh, Le~Thi~Khanh Hien, Nicolas Gillis, and Panagiotis Patrinos.
\newblock A block inertial {B}regman proximal algorithm for nonsmooth nonconvex
  problems with application to nonnegative matrix tri-factorization.
\newblock {\em arXiv:2003.03963}, 2020.

\bibitem{attouch2010proximal}
H{\'e}dy Attouch, J{\'e}r{\^o}me Bolte, Patrick Redont, and Antoine Soubeyran.
\newblock Proximal alternating minimization and projection methods for
  nonconvex problems: An approach based on the {K}urdyka-{{\L}}o\-ja\-sie\-wicz
  inequality.
\newblock {\em Mathematics of Operations Research}, 35(2):438--457, 2010.

\bibitem{attouch2013convergence}
Hedy Attouch, J\'er\^ome Bolte, and Benar~Fux Svaiter.
\newblock Convergence of descent methods for semi-algebraic and tame problems:
  proximal algorithms, forward-backward splitting, and regularized
  {G}auss-{S}eidel methods.
\newblock {\em Mathematical Programming}, 137(1):91--129, Feb 2013.

\bibitem{aze2017nonlinear}
Dominique Az\'e and Jean-No\"el Corvellec.
\newblock Nonlinear error bounds via a change of function.
\newblock {\em Journal of Optimization Theory and Applications}, 172(1):9--32,
  2017.

\bibitem{bauschke2005joint}
Heinz Bauschke, Patrick Combettes, and Dominikus Noll.
\newblock Joint minimization with alternating {B}regman proximity operators.
\newblock {\em Pacific Journal of Optimization}, 2005.

\bibitem{bauschke2019linear}
Heinz~H. Bauschke, J\'er\^ome Bolte, Jiawei Chen, Marc Teboulle, and Xianfu
  Wang.
\newblock On linear convergence of non-{E}uclidean gradient methods without
  strong convexity and {L}ipschitz gradient continuity.
\newblock {\em Journal of Optimization Theory and Applications}, pages 1--20,
  2019.

\bibitem{bauschke2016descent}
Heinz~H. Bauschke, J\'er\^ome Bolte, and Marc Teboulle.
\newblock A descent lemma beyond {L}ipschitz gradient continuity: first-order
  methods revisited and applications.
\newblock {\em Mathematics of Operations Research}, 42(2):330--348, 2016.

\bibitem{bauschke2001essential}
Heinz~H. Bauschke, Jonathan~M. Borwein, and Patrick~L. Combettes.
\newblock Essential smoothness, essential strict convexity, and {L}egendre
  functions in {B}anach spaces.
\newblock {\em Communications in Contemporary Mathematics}, 03(04):615--647,
  2001.

\bibitem{bauschke2018regularizing}
Heinz~H. Bauschke, Minh~N. Dao, and Scott~B. Lindstrom.
\newblock Regularizing with {B}regman-{M}oreau envelopes.
\newblock {\em SIAM Journal on Optimization}, 28(4):3208--3228, 2018.

\bibitem{beck2009fast}
Amir Beck and Marc Teboulle.
\newblock A fast iterative shrinkage-thresholding algorithm for linear inverse
  problems.
\newblock {\em SIAM Journal on Imaging Sciences}, 2(1):183--202, 2009.

\bibitem{bhojanapalli2016global}
Srinadh Bhojanapalli, Behnam Neyshabur, and Nati Srebro.
\newblock Global optimality of local search for low rank matrix recovery.
\newblock In {\em Advances in Neural Information Processing Systems}, pages
  3873--3881, 2016.

\bibitem{bolte2007lojasiewicz}
J\'er\^ome Bolte, Aris Daniilidis, and Adrian Lewis.
\newblock The {{\L}}o\-ja\-sie\-wicz inequality for nonsmooth subanalytic
  functions with applications to subgradient dynamical systems.
\newblock {\em SIAM Journal on Optimization}, 17(4):1205--1223, 2007.

\bibitem{bolte2007clarke}
J\'er\^ome Bolte, Aris Daniilidis, Adrian Lewis, and Masahiro Shiota.
\newblock Clarke subgradients of stratifiable functions.
\newblock {\em SIAM Journal on Optimization}, 18(2):556--572, 2007.

\bibitem{bolte2010characterizations}
J\'er\^ome Bolte, Aris Daniilidis, Olivier Ley, and Laurent Mazet.
\newblock Characterizations of {{\L}}o\-ja\-sie\-wicz inequalities: subgradient
  flows, talweg, convexity.
\newblock {\em Transactions of the American Mathematical Society},
  362(6):3319--3363, 2010.

\bibitem{bolte2017error}
J\'er\^ome Bolte, Trong~Phong Nguyen, Juan Peypouquet, and Bruce~W. Suter.
\newblock From error bounds to the complexity of first-order descent methods
  for convex functions.
\newblock {\em Mathematical Programming}, 165(2):471--507, Oct 2017.

\bibitem{bolte2014proximal}
J{\'e}r{\^o}me Bolte, Shoham Sabach, and Marc Teboulle.
\newblock Proximal alternating linearized minimization for nonconvex and
  nonsmooth problems.
\newblock {\em Mathematical Programming}, 146(1--2):459--494, 2014.

\bibitem{bolte2018first}
J\'er\^ome Bolte, Shoham Sabach, Marc Teboulle, and Yakov Vaisbourd.
\newblock First order methods beyond convexity and {L}ipschitz gradient
  continuity with applications to quadratic inverse problems.
\newblock {\em SIAM Journal on Optimization}, 28(3):2131--2151, 2018.

\bibitem{bot2016inertialTseng}
Radu~Ioan Bo{\c t} and Ern\"o~Robert Csetnek.
\newblock An inertial {T}seng's type proximal algorithm for nonsmooth and
  nonconvex optimization problems.
\newblock {\em Journal of Optimization Theory and Applications},
  171(2):600--616, 2016.

\bibitem{bot2016inertial}
Radu~Ioan Bo{\c t}, Ern{\"o}~Robert Csetnek, and Szil{\'a}rd~Csaba
  L{\'a}szl{\'o}.
\newblock An inertial forward-backward algorithm for the minimization of the
  sum of two nonconvex functions.
\newblock {\em EURO Journal on Computational Optimization}, 4(1):3--25, 2016.

\bibitem{bot2019proximal}
Radu~Ioan Bot, Ernoö~Robert Csetnek, and Dang-Khoa Nguyen.
\newblock A proximal minimization algorithm for structured nonconvex and
  nonsmooth problems.
\newblock {\em SIAM Journal on Optimization}, 29(2):1300--1328, 2019.

\bibitem{bregman1967relaxation}
Lev~M. Bregman.
\newblock The relaxation method of finding the common point of convex sets and
  its application to the solution of problems in convex programming.
\newblock {\em USSR computational mathematics and mathematical physics},
  7(3):200--217, 1967.

\bibitem{davis2018stochastic}
Damek Davis, Dmitriy Drusvyatskiy, and Kellie~J MacPhee.
\newblock Stochastic model-based minimization under high-order growth.
\newblock {\em arXiv preprint arXiv:1807.00255}, 2018.

\bibitem{dennis1974characterization}
John E.~Jr. Dennis and Jorge~J. Mor{\'e}.
\newblock A characterization of superlinear convergence and its application to
  quasi-{N}ewton methods.
\newblock {\em Mathematics of Computation}, 28(126):549--560, 1974.

\bibitem{dennis1977quasi}
John E.~Jr. Dennis and Jorge~J. Mor\'e.
\newblock Quasi-{N}ewton methods, motivation and theory.
\newblock {\em SIAM Review}, 19(1):46--89, 1977.

\bibitem{dontchev2012generalizations}
Asen Dontchev.
\newblock Generalizations of the {D}ennis-{M}oré theorem.
\newblock {\em SIAM Journal on Optimization}, 22(3):821--830, 2012.

\bibitem{dragomir2019fast}
Radu-Alexandru Dragomir, Jérôme Bolte, and Alexandre d'Aspremont.
\newblock Fast gradient methods for symmetric nonnegative matrix factorization.
\newblock {\em arXiv preprint arXiv:1901.10791}, 2019.

\bibitem{dragomir2019quartic}
Radu-Alexandru Dragomir, Alexandre d'Aspremont, and J\'er\^ome Bolte.
\newblock Quartic first-order methods for low rank minimization.
\newblock {\em arXiv preprint arXiv:1901.10791}, 2019.

\bibitem{facchinei2003finite}
Francisco Facchinei and Jong-Shi Pang.
\newblock {\em Finite-dimensional variational inequalities and complementarity
  problems}, volume~II.
\newblock Springer, 2003.

\bibitem{fischer2002local}
Andreas Fischer.
\newblock Local behavior of an iterative framework for generalized equations
  with nonisolated solutions.
\newblock {\em Mathematical Programming}, 94(1):91--124, 2002.

\bibitem{fukushima1981generalized}
Masao Fukushima and Hisashi Mine.
\newblock A generalized proximal point algorithm for certain non-convex
  minimization problems.
\newblock {\em International Journal of Systems Science}, 12(8):989--1000,
  1981.

\bibitem{fukushima1996globally}
Masao Fukushima and Liqun Qi.
\newblock A globally and superlinearly convergent algorithm for nonsmooth
  convex minimization.
\newblock {\em SIAM Journal on Optimization}, 6(4):1106--1120, 1996.

\bibitem{hanzely2018fastest}
Filip Hanzely and Peter Richt\'arik.
\newblock Fastest rates for stochastic mirror descent methods.
\newblock {\em arXiv preprint arXiv:1803.07374}, 2018.

\bibitem{hanzely2018accelerated}
Filip Hanzely, Peter Richtarik, and Lin Xiao.
\newblock Accelerated {B}regman proximal gradient methods for relatively smooth
  convex optimization.
\newblock {\em arXiv preprint arXiv:1808.03045}, 2018.

\bibitem{izmailov2014newton}
Alexey~F. Izmailov and Mikhail~V. Solodov.
\newblock {\em {N}ewton-type {M}ethods for {O}ptimization and {V}ariational
  {P}roblems}.
\newblock Springer, 2014.

\bibitem{kan2012moreau}
Chao Kan and Wen Song.
\newblock The {M}oreau envelope function and proximal mapping in the sense of
  the {B}regman distance.
\newblock {\em Nonlinear Analysis: Theory, Methods \& Applications}, 75(3):1385
  -- 1399, 2012.

\bibitem{kawaguchi2016deep}
Kenji Kawaguchi.
\newblock Deep learning without poor local minima.
\newblock In {\em Advances in neural information processing systems}, pages
  586--594, 2016.

\bibitem{kurdyka1998gradients}
Krzysztof Kurdyka.
\newblock On gradients of functions definable in o-minimal structures.
\newblock {\em Annales de l'institut Fourier}, 48(3):769--783, 1998.

\bibitem{laude2020bregman}
Emanuel Laude, Peter Ochs, and Daniel Cremers.
\newblock {B}regman proximal mappings and {B}regman--{M}oreau envelopes under
  relative prox-regularity.
\newblock {\em Journal of Optimization Theory and Applications},
  184(3):724--761, 2020.

\bibitem{lee2019first}
Jason~D. Lee, Ioannis Panageas, Georgios Piliouras, Max Simchowitz, Michael~I.
  Jordan, and Benjamin Recht.
\newblock First-order methods almost always avoid strict saddle points.
\newblock {\em Mathematical Programming}, Feb 2019.

\bibitem{li2004regularized}
Dong-Hui Li, Masao Fukushima, Liqun Qi, and Nobuo Yamashita.
\newblock Regularized {N}ewton methods for convex minimization problems with
  singular solutions.
\newblock {\em Computational Optimization and Applications}, 28(2):131--147,
  2004.

\bibitem{li2016douglas}
Guoyin Li and Ting~Kei Pong.
\newblock Douglas-{R}achford splitting for nonconvex optimization with
  application to nonconvex feasibility problems.
\newblock {\em Mathematical Programming}, 159(1):371--401, Sep 2016.

\bibitem{li2018highly}
Xudong Li, Defeng Sun, and Kim-Chuan Toh.
\newblock A highly efficient semismooth {N}ewton augmented {L}agrangian method
  for solving lasso problems.
\newblock {\em SIAM Journal on Optimization}, 28(1):433--458, 2018.

\bibitem{liu2017further}
Tianxiang Liu and Ting~Kei Pong.
\newblock Further properties of the forward-backward envelope with applications
  to difference-of-convex programming.
\newblock {\em Computational Optimization and Applications}, 67(3):489--520,
  Jul 2017.

\bibitem{liu2019envelope}
Yanli Liu and Wotao Yin.
\newblock An envelope for {D}avis--{Y}in splitting and strict saddle-point
  avoidance.
\newblock {\em Journal of Optimization Theory and Applications},
  181(2):567--587, May 2019.

\bibitem{lojasiewicz1963propriete}
Stanislaw {\L}ojasiewicz.
\newblock Une propri{\'e}t{\'e} topologique des sous-ensembles analytiques
  r{\'e}els.
\newblock {\em Les {\'e}quations aux d{\'e}riv{\'e}es partielles}, pages
  87--89, 1963.

\bibitem{lojasiewicz1993geometrie}
Stanislaw {\L}ojasiewicz.
\newblock Sur la g{\'e}om{\'e}trie semi- et sous- analytique.
\newblock {\em Annales de l'institut Fourier}, 43(5):1575--1595, 1993.

\bibitem{lu2018relatively}
Haihao Lu, Robert~M. Freund, and Yurii Nesterov.
\newblock Relatively smooth convex optimization by first-order methods, and
  applications.
\newblock {\em SIAM Journal on Optimization}, 28(1):333--354, 2018.

\bibitem{maratos1978exact}
Nicholas Maratos.
\newblock {\em Exact penalty function algorithms for finite dimensional and
  control optimization problems}.
\newblock PhD thesis, Imperial College London (University of London), 1978.

\bibitem{mukkamala2019beyond}
Mahesh~Chandra Mukkamala and Peter Ochs.
\newblock Beyond alternating updates for matrix factorization with inertial
  {B}regman proximal gradient algorithms.
\newblock In {\em Advances in Neural Information Processing Systems}, pages
  4268--4278, 2019.

\bibitem{mukkamala2019convex}
Mahesh~Chandra Mukkamala, Peter Ochs, Thomas Pock, and Shoham Sabach.
\newblock Convex-concave backtracking for inertial {B}regman proximal gradient
  algorithms in non-convex optimization.
\newblock {\em arXiv preprint arXiv:1904.03537}, 2019.

\bibitem{nesterov2013gradient}
Yurii Nesterov.
\newblock Gradient methods for minimizing composite functions.
\newblock {\em Mathematical Programming}, 140(1):125--161, Aug 2013.

\bibitem{nesterov2018implementable}
Yurii Nesterov.
\newblock Implementable tensor methods in unconstrained convex optimization.
\newblock Technical report, UC Louvain, Center for Operations Research and
  Econometrics (CORE), Belgium, 2018.

\bibitem{nocedal2006numerical}
Jorge Nocedal and Stephen Wright.
\newblock {\em Numerical {O}ptimization}.
\newblock Springer Science \& Business Media, 2006.

\bibitem{noll2014convergence}
Dominikus Noll.
\newblock Convergence of non-smooth descent methods using the
  kurdyka--{\l}ojasiewicz inequality.
\newblock {\em Journal of Optimization Theory and Applications},
  160(2):553--572, 2014.

\bibitem{ochs2019nonsmooth}
Peter Ochs, Jalal Fadili, and Thomas Brox.
\newblock Non-smooth non-convex {B}regman minimization: Unification and new
  algorithms.
\newblock {\em Journal of Optimization Theory and Applications},
  181(1):244--278, 2019.

\bibitem{oneill2018behavior}
Michael O'Neill and Stephen~J. Wright.
\newblock Behavior of accelerated gradient methods near critical points of
  nonconvex functions.
\newblock {\em Mathematical Programming}, Oct 2018.

\bibitem{palis1982geometric}
Jacob~Jr Palis and Welington De~Melo.
\newblock {\em Geometric theory of dynamical systems: an introduction}.
\newblock Springer-Verlag, 1982.

\bibitem{patrinos2013proximal}
Panagiotis Patrinos and Alberto Bemporad.
\newblock Proximal {N}ewton methods for convex composite optimization.
\newblock In {\em 52nd IEEE Conference on Decision and Control}, pages
  2358--2363, 2013.

\bibitem{poliquin1995second}
Ren{\'e}~A. Poliquin and R.~Tyrrell Rockafellar.
\newblock Second-order nonsmooth analysis in nonlinear programming.
\newblock {\em Recent advances in nonsmooth optimization}, pages 322--349,
  1995.

\bibitem{poliquin1996proxregular}
Ren\'e~A. Poliquin and R.~Tyrrell Rockafellar.
\newblock Prox-regular functions in variational analysis.
\newblock {\em Transactions of the American Mathematical Society},
  348(5):1805--1838, 1996.

\bibitem{rockafellar1977higher}
R.~Tyrrell Rockafellar.
\newblock Higher derivatives of conjugate convex functions.
\newblock {\em Journal of Applied Analysis}, 1(1):41--43, 1977.

\bibitem{rockafellar2011variational}
R.~Tyrrell Rockafellar and Roger J.-B. Wets.
\newblock {\em Variational {A}nalysis}, volume 317.
\newblock Springer Science \& Business Media, 2011.

\bibitem{rockafellar1970convex}
Ralph~Tyrell Rockafellar.
\newblock {\em Convex {A}nalysis}.
\newblock Princeton university press, 1970.

\bibitem{stella2017forward}
Lorenzo Stella, Andreas Themelis, and Panagiotis Patrinos.
\newblock Forward-backward quasi-{N}ewton methods for nonsmooth optimization
  problems.
\newblock {\em Computational Optimization and Applications}, 67(3):443--487,
  Jul 2017.

\bibitem{stella2017simple}
Lorenzo Stella, Andreas Themelis, Pantelis Sopasakis, and Panagiotis Patrinos.
\newblock A simple and efficient algorithm for nonlinear model predictive
  control.
\newblock In {\em 2017 IEEE 56th Annual Conference on Decision and Control
  (CDC)}, pages 1939--1944, Dec 2017.

\bibitem{sun2018geometric}
Ju~Sun, Qing Qu, and John Wright.
\newblock A geometric analysis of phase retrieval.
\newblock {\em Foundations of Computational Mathematics}, 18(5):1131--1198,
  2018.

\bibitem{tanner2016low}
Jared Tanner and Ke~Wei.
\newblock Low rank matrix completion by alternating steepest descent methods.
\newblock {\em Applied and Computational Harmonic Analysis}, 40(2):417--429,
  2016.

\bibitem{teboulle1992entropic}
Marc Teboulle.
\newblock Entropic proximal mappings with applications to nonlinear
  programming.
\newblock {\em Mathematics of Operations Research}, 17(3):670--690, 1992.

\bibitem{teboulle2018simplified}
Marc Teboulle.
\newblock A simplified view of first order methods for optimization.
\newblock {\em Mathematical Programming}, pages 1--30, 2018.

\bibitem{themelis2018proximal}
Andreas Themelis.
\newblock {\em Proximal Algorithms for Structured Nonconvex Optimization}.
\newblock PhD thesis, KU Leuven / IMT Lucca, 2018.

\bibitem{themelis2019acceleration}
Andreas Themelis, Masoud Ahookhosh, and Panagiotis Patrinos.
\newblock On the acceleration of forward-backward splitting via an inexact
  {N}ewton method.
\newblock In R.~Luke, H.~Bauschke, and R.~Burachik, editors, {\em Splitting
  Algorithms, Modern Operator Theory, and Applications}, pages 363--412.
  Springer, 2019.

\bibitem{themelis2019supermann}
Andreas Themelis and Panagiotis Patrinos.
\newblock {S}uper{M}ann: a superlinearly convergent algorithm for finding fixed
  points of nonexpansive operators.
\newblock {\em IEEE Transactions on Automatic Control}, dec 2019.

\bibitem{themelis2018forward}
Andreas Themelis, Lorenzo Stella, and Panagiotis Patrinos.
\newblock Forward-backward envelope for the sum of two nonconvex functions:
  Further properties and nonmonotone linesearch algorithms.
\newblock {\em SIAM Journal on Optimization}, 28(3):2274--2303, 2018.

\bibitem{ueda2014regularized}
Kenji Ueda and Nobuo Yamashita.
\newblock A regularized {N}ewton method without line search for unconstrained
  optimization.
\newblock {\em Computational Optimization and Applications}, 59(1-2):321--351,
  2014.

\bibitem{van1998tame}
Lou Van~den Dries.
\newblock {\em Tame topology and o-minimal structures}, volume 248.
\newblock Cambridge university press, 1998.

\bibitem{van1996geometric}
Lou Van~den Dries and Chris Miller.
\newblock Geometric categories and o-minimal structures.
\newblock {\em Duke Mathematical Journal}, 84(2):497--540, 08 1996.

\bibitem{van2017forward}
Quang Van~Nguyen.
\newblock Forward-backward splitting with {B}regman distances.
\newblock {\em Vietnam Journal of Mathematics}, 45(3):519--539, 2017.

\bibitem{yamashita2001rate}
Nobuo Yamashita and Masao Fukushima.
\newblock On the rate of convergence of the {L}evenberg-{M}arquardt method.
\newblock In {\em Topics in numerical analysis}, pages 239--249. Springer,
  2001.

\bibitem{yu2019deducing}
Peiran Yu, Guoyin Li, and Ting~Kei Pong.
\newblock Deducing {K}urdyka-{{\L}}o\-ja\-sie\-wicz exponent via
  inf-projection.
\newblock {\em arXiv preprint arXiv:1902.03635}, 2019.

\bibitem{yue2019quadratic}
Man-Chung Yue, Zirui Zhou, and Anthony Man-Cho~So.
\newblock On the quadratic convergence of the cubic regularization method under
  a local error bound condition.
\newblock {\em SIAM Journal on Optimization}, 29(1):904--932, 2019.

\end{thebibliography}

\end{document}